\documentclass[preprint,sort&compress,12pt]{elsarticle}
\usepackage{bbm}
\usepackage{amsmath}
\usepackage{mathrsfs}
\usepackage{stmaryrd}
\usepackage{mathrsfs}
\usepackage{bm}
\usepackage{amsmath}
\usepackage{amsfonts}
\usepackage{amssymb}
\usepackage{amsthm}
\usepackage{lineno}
\usepackage[RGB]{xcolor}
\newcommand{\mm}{\mathrm}
\newcommand{\ml}{\mathcal}

\newcommand{\be}{\begin{equation}}
\newcommand{\bea}{\begin{equation}\begin{aligned}}
\newcommand{\beas}{\begin{equation*}\begin{aligned}}
\newcommand{\eeas}{\end{aligned}\end{equation*}}
\newcommand{\eea}{\end{aligned}\end{equation}}
\newcommand{\ee}{\end{equation}}

\renewcommand{\div}{{\rm div }}

\begin{document}
\begin{frontmatter}
\title{On Magnetic Inhibition Theory
in Non-resistive Magnetohydrodynamic Fluids: Existence of Solutions in Some Classes of Large Data}
% \cortext[cor1]
%{We  certify that the general content of this article, in whole or in part, is not submitted, accepted, or published elsewhere, including conference proceedings. }

\author[FJ,AB]{Fei Jiang}\ead{jiangfei0591@163.com}
% \cortext[cor1]{We  certify that the general content of this article, in whole or in part, is not submitted, accepted, or published elsewhere, including conference proceedings. }
\author[sJ]{Song Jiang}\ead{jiang@iapcm.ac.cn}
%\author[FJ]{Yanjie Zhang}\ead{zhangyajie315@qq.com}
\address[FJ]{School of Mathematics and Statistics, Fuzhou University, Fuzhou, 350108, China.}
\address[sJ]{Institute of Applied Physics and Computational Mathematics, Huayuan Road 6, Beijing 100088, China.}
\address[AB]{Key Laboratory of Operations Research and Control of Universities in Fujian, Fuzhou 350108, China.}
\begin{abstract}
This paper is concerned with existence of solutions  to the incompressible non-resistive viscous magnetohydrodynamic (MHD) equations with large initial perturbations in there-dimensional (3D) periodic domains (in Lagrangian coordinates). Motivated by the Diophantine condition imposed by the approximate theory of non-resistive MHD equations in \cite{BCSCSPLL},  Chen--Zhang--Zhou in \cite{chen2021} and  the magnetic inhibition mechanism  of Lagrangian coordinates version in our previous paper \cite{JFJSOMITIN}, we prove the existence of unique classical solutions under some class of large initial perturbations, where the intensity of impressive magnetic fields depends increasingly on the $ H^{17}\times H^{21}$-norm of the initial perturbation of both the velocity and magnetic field. Our result not only mathematically  verifies that magnetic fields prevent the singularity formation  of  solutions with large initial velocity in the viscous case, but also provide a starting point for the existence theory of large perturbation solutions of the 3D non-resistive viscous MHD equations. In addition, we further rigorously prove that, for large time or strong magnetic field, the MHD equations reduce to the corresponding linearized equations by providing the error estimates, which enjoy the types of algebraic decay with respect to the both of time and field intensity, between the solutions of both the nonlinear and linear equations.
\end{abstract}
\begin{keyword}
Incompressible MHD fluids;  large initial data; algebraic time-decay;  {convergence in the field intensity}.
\MSC[2000] 35Q60\sep  35B10\sep 76E25.
%(2000 is the default)
\end{keyword}
\end{frontmatter}

%% Start linex numbering here if you wantxx
% \linenumbers

%% main text
\newtheorem{thm}{Theorem}[section]
\newtheorem{lem}{Lemma}[section]
\newtheorem{pro}{Proposition}[section]
\newtheorem{cor}{Corollary}[section]
\newproof{pf}{Proof}
\newdefinition{rem}{Remark}[section]
\newtheorem{definition}{Definition}[section]
% \linenumbers
\newtheorem{concl}{Assertion}[section]
\section{Introduction}\label{introud}
\numberwithin{equation}{section}

In this paper, we investigate the existence of the global (-in-time) solutions to the following system of equations of  {an} incompressible non-resistive viscous
magnetohydrodynamic (MHD) fluid  with large initial perturbations:
\begin{equation}\label{1.1xx}
\begin{cases}
\rho  {v}_t+\rho {v}\cdot\nabla {v}+\nabla p- \mu\Delta v=\lambda M \cdot\nabla M/4\pi, \\
M_t+v\cdot\nabla M =M\cdot \nabla v,\\
\mathrm{div} v=\mathrm{div} M=0,
\end{cases}
\end{equation}
where the unknowns ${v}:={v}(x,t)$, $M:=M(x,t)$ and $p:=p(x,t)$ denote the velocity, magnetic field and the sum of both  {the}
magnetic and
kinetic pressures of  {the MHD fluid}, resp., and the three positive (physical) parameters $\rho$, $\mu$ and $\lambda$
stand for the density, shear viscosity coefficient and permeability of vacuum, resp.. In particular, neglecting the shear viscosity coefficient, we have
\begin{equation}\label{1x}
\begin{cases}
\rho  {v}_t+\rho {v}\cdot\nabla {v}+\nabla p =\lambda M \cdot\nabla M/4\pi, \\
M_t+v\cdot\nabla M =M\cdot \nabla v,\\
\mathrm{div} v=\mathrm{div} M=0.
\end{cases}
\end{equation}
%Exploiting the Els\"asser variables $$Z^+:=v+M,\ Z^-:=v-M,\ $$

Physicists pointed out that, in the nonlinear MHD system, a strong enough magnetic field will reduce the nonlinear interaction \cite{kraichnan1965inertial} and inhibit the formation of strong gradients. This effect was also observed in direct numerical simulations of the ideal MHD system \eqref{1x}, with periodic boundary conditions \cite{frisch1983dynamics}.
In 1998, Bardos--Sulem--Sulem used the hyperbolicity of the two/three-dimensional (2/3D) idea MHD system \eqref{1x} and gave a rigorous proof of the global well-posedness when the initial data $(v^0,M^0)$ is around the equilibrium state $(0,\bar{M})$ in  the H\"older space \cite{BCSCSPLL} (see \cite{CYLZGW,HLBXLYP,xuli2020,Dongyi2017} for the case of Sobolev spaces), where $\bar{M}$ is often called the impressive magnetic field. In particular, the result of Bardos et.al. presents that the system \eqref{1x} is global well-posed with large initial perturbation, when the impressive magnetic field is sufficiently large. We remark that such global stability result is not excepted for the 3D  incompressible Euler equations, i.e., the magnetic field is absent in the system \eqref{1x}. Interesting readers can refer to   \cite{luo2014potentially,CGENJJGTY2021} and \cite{kiselev2014small,zlatovs2015exponential} for the  singularity formation of solutions  and the solutions with double-exponential growth in Euler equations, resp.. Bardos--Sulem--Sulem's large perturbation result can  be roughly described as follows:
\begin{concl}If
 \begin{align}
 \frac{\mbox{initial perturbation around $(0,\bar{M})$}}{\mbox{field intensity of $\bar{M}$}}\ll 1,
\nonumber
 \end{align}
 then the idea incompressible MHD system \eqref{1x} admits a global stability solution. Here and in what follows, $A\ll B$ means that $A$ is much smaller than $B$.
\end{concl}
We naturally make  an association with the following well-known assertion in viscous (pure) fluids (see  \cite{lei2011global,fujita1964navier,craig2013global} for examples):
\begin{concl}If
 \begin{align}
 \frac{\mbox{initial perturbation of velocity around zero}}{\mbox{value of viscocity } \mu}\ll 1,\nonumber
 \end{align}
 then Navier--Stokes equations (i.e. the equations in \eqref{1.1xx} with $M=0$) admits a global stability solution.
\end{concl}

The above two assertions present  that magnetic fields can inhibit the singularity formation of solutions  with large initial velocity as well as viscosity, though the physical mechanisms of inhibition/stabilzing are different. Based on the above two assertions, we easily further believe that such large perturbation result shall also exists for the hyperbolic-parabolic system \eqref{1.1xx}, that is (roughly speaking):
 \begin{concl}
 \label{2020111031231}
 If
 \begin{align}
 \frac{\mbox{initial perturbation around $(0,\bar{M})$}}{\mbox{field intensity of $\bar{M}$}}\ll 1,
 \label{20202111011816}
 \end{align}
 then the  incompressible non-resistive viscous MHD system \eqref{1.1xx} admits a global stability solution.
\end{concl}

Unfortunately, due to the parabolic structure of \eqref{1.1xx}$_1$, the hyperbolic method used by Bardos et.al. can not be applied to the system \eqref{1.1xx}. Nearly twenty years after the pioneering work of Bardos et.al., Zhang developed new ideas based on the energy method and spectrum analysis, and successfully established the large perturbation result for the system \eqref{1.1xx} under the 2D case \cite{ZTGS}. Zhang's result supports Assertion \ref{2020111031231}, at lest for the 2D case. {However the corresponding 3D case is still an open problem}. At present, all 3D global well-posedness results without symmetry structure imposed on the initial data (see \cite{lei2015axially} for the axially symmetric solutions with large initial data) are concerned with the small perturbation solutions, i.e.  if \begin{align}
 {\mbox{initial perturbation around $(0,\bar{M})$}} \ll 1,
 \end{align}
  then the system \eqref{1.1xx} admits a global stability solution,  see \cite{xu2015global,ABIHZPOTG,pan2018global} and \cite{TZWYJGw} for the Cauchy problem and the initial-boundary value
 problem of \eqref{1.1xx}, resp.. Interesting readers can further refer to \cite{RXXWJHXZYZZF,lin2015global,zhang2014elementary,ren2016global} and \cite{fefferman2017local,HOCELETFCMDSFA,chemin2016local} for the 2D small perturbation results and the 3D local (-in-time) existence results of large initial perturbations, resp.. In this paper, we will develop some new ideas to further establish the existence of large perturbation solutions for the system \eqref{1.1xx} in Lagrangian coordinates under a relatively small condition as \eqref{20202111011816}, see Remark \ref{2021110321826} for further discussion.  Our result  provides a starting point for the existence theory of large perturbation solutions of the 3D non-resistive viscous MHD fluids. Next let us first recall a heuristic physical idea, which lights up the road to  Assertion \ref{2020111031231} for the 3D case.

 It is well-known that we can think that the 3D non-resistive MHD fluid under rest state is made up of infinite (fluid) element lines which are parallel to the impressive field and can be regarded  elastic strings under magnetic tension. Let us pick up a segment from an element line, and denote it by $y^1y^2$, see the following first figure, in which $F_{y^i}$ denotes the magnetic tension acting on the element endpoint $y^i$ for $i=1$, $2$, and is given by the formula (see  \cite{JFJSOMITIN} for the detailed derivation)
\begin{align}
\label{20211021009}
F_{y^i}=\lambda  \partial_{\bar{M}} \zeta(y^i,s)\cdot \nu_{y^i}   \partial_{\bar{M}} \zeta(y^i,s)/4\pi,
\end{align}
 where $\zeta$ denotes the flow function (see \eqref{201806122101} for the definition) and $\nu_{y^i}$ represents the unit outer normal vector at element endpoint $y^i$. Once we disturb the rest state, the segment $y^1y^2$  will be bend, see the following second figure, in which $F_{\mm{r}}$ denotes the resultant force of $F_{y^1}$ and $F_{y^2}$.  However the magnetic tension will  straighten the  bent element line, playing the role of resilience. Since
   the magnetic tension intensity is strictly increasing on the impressed field intensity by the formula \eqref{20211021009}, we easily see that a strong enough magnetic field can inhibit the singularity formation, and even the flow instabilities \cite{JFJSWWWOA,JFJSJMFMOSERT,JFJSSETEFP,JFJSJMFM,JFJSOUI,WYTIVNMI,CSHHSCPO}. Moreover, such inhibition mechanism is also  useful for understanding the stabilizing effect of magnetic fields on the motion of resistive MHD fluids \cite{GGPNA,ZYZYGC,WYJCMP}
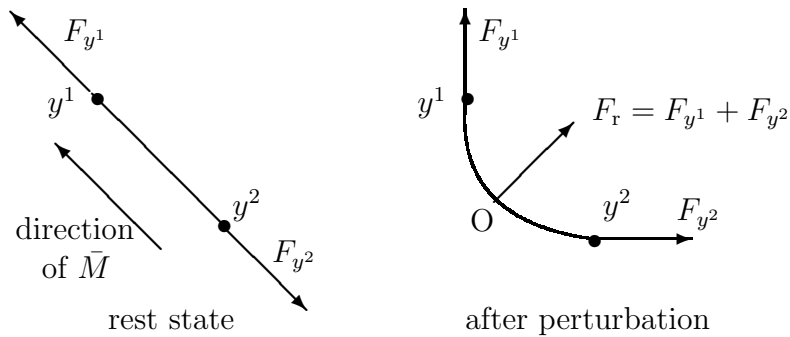
\begin{center}
\begin{picture}(0,130)(128,0)
\thicklines
\put(28,84){$\bullet$}
\put(12,82){$y^1$}
\put(18,109){$F_{y^1}$}
\put(28,90){\vector(-1,1){30}}
\put(76,36){$\bullet$}
\put(29,89){\line(1,-1){50}}
\put(82,43){$y^2$}
\put(97,25){$F_{y^2}$}
\put(80,37){\vector(1,-1){30}}
\put(55,30){\vector(-1,1){40}}
\put(0,34){direction}
\put(10,19){of $\bar{M}$}
\put(35,0){rest state}
\put(168,84){$\bullet$}
\put(152,82){$y^1$}
\put(170,81){\vector(0,1){40}}
\put(175,109){$F_{y^1}$}
\put(216,30){$\bullet$}
\put(222,43){$y^2$}
\put(216,34){\vector(1,0){40}}
\put(250,41){$F_{y^2}$}
\put(172,37){O}
\put(181,48){\vector(1,1){30}}
\put(218,79){$F_{\mm{r}}= F_{y^1}+F_{y^2}$}
\qbezier(170,84)(167,41)(218,34)
\put(170,0){after perturbation}
\end{picture}
\end{center}

Since the above inhibition mechanism is described in Lagrangian coordinates, this naturally motivates us to consider the large perturbation condition \eqref{20202111011816} in Lagrangian coordinates. Hence we write down the well-known  energy law of \eqref{1.1xx} in Lagrangian coordinates:
 \begin{align}
& \rho\int_\Omega |u|^2\mm{d}y+ \frac{\lambda \varpi^2}{4\pi}\int_\Omega |\partial_{\omega}\eta|^2\mm{d}y
+ {2\mu} \int_0^t \int_\Omega |\nabla_{\mathcal{A}} u |^2 \mm{d}y \mm{d}\tau\nonumber \\
&= \rho\int_\Omega |u^0|^2 \mm{d}y + \frac{\lambda \varpi^2}{4\pi}\int_\Omega |\partial_{ \omega} \eta^0 |^2 \mm{d}y =:I^0.  \label{202211121103}
\end{align}

We shall explain the notations appearing in the above identity.
 $\Omega$ is the fluid domain.
 We have defined that \begin{align}
\bar{M}:=\varpi\omega,
\label{202111031455}
 \end{align}where $\omega$ and $m$ denote the unit vector in $\mathbb{R}^3$ and the field intensity of $\bar{M}$, resp..  $\eta:=\zeta-y$ represents the deviation function of particles, and $u$ is the velocity function in Lagrangian coordinates. The differential operator $\nabla_{\mathcal{A}}$ is given by \eqref{202111021116}. Here and in what follows, $f^0$ always denotes the initial data of $f$. We call $I^0$ the initial mechanical energy \cite{JFJSOMITIN}.

It is easy to see from \eqref{202211121103} that
  \begin{align}
 {\partial_{\omega}\eta} \ll 1,
\end{align}
if the initial mechanical energy $I_0$ is given and satisfies
\begin{align}
\label{2021110221127}
\frac{\sqrt{I^0}}{\varpi }\ll 1.
\end{align}
The above fact provides a new idea to prove the existence of (global) large perturbation solutions for the system \eqref{1.1xx}. Motivated by \eqref{2021110221127} and the odevity conditions imposed by Pan--Zhou--Zhu in \cite{pan2018global}, for the 2D spatially periodic domain $\mathbb{T}^2$ and $\omega=(0,1)^{\mm{T}}$, the authors renew to prove the the existence of large perturbation solutions for the system \eqref{1.1xx} in Lagrangian coordinates \cite{JFJSAAD}, where $\mathbb{T}:=\mathbb{R}\setminus\mathbb{Z}$ and the initial data satisfies some relatively small condition similar to \eqref{20202111011816}.

The significance of odevity condition lies in that the following important inverse  relation can be further established for the 2D case:
  \begin{align}
  \label{202210}
  {\nabla \eta} \propto  \varpi ^{-1} .
\end{align}
This relation intuitively reveals that the (nonlinear) solutions of \eqref{1.1xx} in Lagrangian coordinates can be approximated
by the (linear) solutions of the corresponding linearized equations for $\varpi \gg 1$, and thus we can expect to establish the existence of large perturbation solutions. However it seems  that this idea by using odevity conditions fails to the corresponding 3D case, since the relation  \eqref{202210} for the 3D case can not be obtained.

Recently Chen--Zhang--Zhou made an important progress \cite{chen2021}. They observed that if $\bar{M}$  satisfies the Diophantine condition:
\begin{align}
\mbox{there exists a constant }  c_{\bar{M}}>0\mbox{ such that }  |\chi\cdot\bar{M}|\geqslant c_{\bar{M}}|\chi|^{-3}\mbox{ for any }\chi\in \mathbb{Z}^3\setminus \{0\},
\label{202755}
\end{align}
where   $c_{\bar{M}}$ depends on $\bar{M}$ and $\cdot$ denotes the inner product of two vectors,
then one has the generalized  Poinc\'are's inequality \eqref{2020211013200910}, which results into that
the existence of small perturbation solutions for the system \eqref{1.1xx}  can be directly proved without the help of any transformation of functions. Motivated by the Chen--Zhang--Zhou's work, we amazingly find that if $\bar{M}=\varpi \omega$ and the unit vector $\omega\in \mathbb{R}^3$ satisfies  the Diophantine condition (please refer to Remark \ref{202110121947} for the existence)
\begin{align}
\exists\mbox{ a constnat }c_{\omega}>0,\mbox{ s.t. }|\chi\cdot\omega|\geqslant c_{\omega}|\chi|^{-3}\mbox{ for any }\chi\in \mathbb{Z}^3\setminus \{0\},
\label{20201102121755}
\end{align}
then the relation \eqref{202210} also holds for the 3D case, and thus the existence of large perturbation solutions for the system \eqref{1.1xx} can be establish by directly using a multi-layers energy method, see Theorem \ref{201904301948} for the details.

We mention that the proof of our 3D result in Theorem \ref{201904301948} is very different to Zhang's 2D one in \cite{ZTGS}.  In fact, Zhang first obtained a (linear)  solution of linearized equations in Eulerian coordinates, then proved the existence of {a small} error solution between the both linear and nonlinear solutions for strong magnetic field, and finally got a large solution   by adding the linear solution and the small error solution together.
 However the relation \eqref{202210} for the 3D case allows us to directly establish the existence of solutions
  under some class of large initial perturbations by one-step procedure in this paper, rather than Zhang's three-step procedure.
   In addition,  we further rigorously prove that, for large time or strong magnetic field, the MHD equations reduce to the corresponding linearized equations. More precisely,
\begin{itemize}
\item the  nonlinear interactions of the large perturbation solutions become asymptotically negligible as $t\to \infty$ in Lagrangian coordinates  (Such phenomenon had been verified by Bardos--Sulem--Sulem for the system \eqref{1x} \cite{BCSCSPLL});
\item  the difference between the both solutions of \eqref{1.1xx} in Lagrangian coordinates and the corresponding linearized system
can be bounded from above by $ O(\varpi ^{-1/2})$ as $\varpi \to \infty$.
\end{itemize}
please refer to Theorem \ref{201912041028} for details.
Finally, we mention that the asymptotic behaviors of solutions with respect to other parameters, such as the Mach and Alfv\'en numbers,
in MHD fluids have been also extensively investigated, see, for example, \cite{CBHQCSS} and the references cited therein.

The rest of this paper is organized as follows: In Section \ref{202008161246} we introduce our main results including the existence of unique classical solutions with some class of large initial data to the 3D system \eqref{1.1xx} in a periodic domain in Lagrangian coordinates and the both convergence rates of the classical solutions for $\varpi \to\infty$ and $t\to\infty$, i.e., Theorems \ref{201904301948}  and \ref{201912041028}, the proofs of which are given in Sections \ref{20190sdafs4301948}--\ref{2020011sdf92326} in sequence.

\section{Main results}\label{202008161246}
In this section we will state the main results in details. To begin with, we reformulate the equations \eqref{1.1xx} in Lagrangian coordinates.  {Recalling that we study \eqref{1.1xx} in a 3D periodic domain in this paper}, we see, without loss of generality,
that it suffices to consider the periodic domain $\mathbb{T}^3$ with $\mathbb{T}:= \mathbb{R}/\mathbb{Z}$.
\subsection{Reformulation in  Lagrangian coordinates}
Let $(v,M)$ be the solution of  {the 3D system \eqref{1.1xx}}, and the flow map $\zeta$ be the solution to
\begin{equation}
\label{201806122101}
            \begin{cases}
\partial_t \zeta(y,t)=v(\zeta(y,t),t)&\mbox{in }\mathbb{T}^3 \times\mathbb{R}^+ ,
\\
\zeta(y,0)=\zeta^0(y)&\mbox{in }\mathbb{T}^3,
                  \end{cases}
\end{equation}
where $\mathbb{R}^+=(0,\infty)$, $\zeta^0(y)$ satisfies $\det \nabla \zeta^0=1$ and ``$\det$" denotes the determinant.

Since $v$ is divergence-free, then
\begin{align}
\det \nabla \zeta =1
\label{202008161709}
\end{align}
as well as $\det \nabla \zeta^0=1$. Thus, we define
\begin{align}\ml{A}^{\mm{T}}:=(\nabla \zeta)^{-1}:=(\partial_j \zeta_i)^{-1}_{3\times 3},\nonumber
\end{align}
where the superscript $\mm{T}$ represents transposition.

 {Now, we introduce some} differential operators involving $\mathcal{A}$, which will be used later. The differential
 operators $\nabla_{\ml{A}}$, $\mm{div}_\ml{A}$ and $\Delta_{\ml{A}}$ are defined by
 \begin{align}
 & \nabla_{\ml{A}}f:=(\ml{A}_{1k}\partial_kf,\ml{A}_{2k}\partial_kf ,
 \ml{A}_{3k}\partial_kf)^{\mm{T}},
\label{202111021116} \\
& \mm{div}_{\ml{A}}(X_1,X_2,X_3)^{\mm{T}}:=\ml{A}_{lk}\partial_k X_l\mbox{ and } \Delta_{\mathcal{A}}f:=\mm{div}_{\ml{A}}\nabla_{\ml{A}}f\nonumber
\end{align}
for a scalar function $f$ and a vector function $X:=(X_1,X_2,X_3)^{\mm{T}}$, where $\ml{A}_{ij}$ denotes the $(i,j)$-th entry of the matrix $\mathcal{A}$. It should be remarked that  we have used the Einstein convention of summation over repeated indices, and $\partial_k=\partial_{y_k}$.
In addition, thanks to \eqref{202008161709}, we have
\begin{equation}
\partial_k  \mathcal{A}_{ik} =0.
\label{201909261909}
\end{equation}

Let  $\nu= \mu/{\rho}$ and
\begin{equation*}
( u ,B,q)(y,t)=( v,M,p/\rho)(\zeta(y,t),t) \mbox{ for }  (y,t)\in \mathbb{T}^3 \times\mathbb{R}^+ .
\end{equation*}
By virtue of the equations \eqref{1.1xx} and \eqref{201806122101}$_1$, the evolution equations for $(\zeta,u,q)$ in Lagrangian
coordinates read as follows.
\begin{equation}\label{01dsaf16asdfasfdsaasf}
                              \begin{cases}
\zeta_t=u , \\[1mm]
 u_t+\nabla_{\ml{A}}q- \nu \Delta_{\ml{A}}u=
 \lambda B\cdot\nabla_{\mathcal{A}}B/4\pi\rho, \\[1mm]
B_t-B\cdot\nabla_{\mathcal{A}} u=0,\\
\div_\ml{A}u=0,\\[1mm]\div_\ml{A}B=0.
\end{cases}
\end{equation}

 We can  derive from \eqref{01dsaf16asdfasfdsaasf}$_3$ the differential version of magnetic flux conservation \cite{JFJSOMITIN}:
\begin{equation} \nonumber
\ml{A}_{jl}B_j=\ml{A}_{jl}^0B_j^0,
 \end{equation}
 which yields
\begin{eqnarray}
 \label{0124}  B=\nabla\zeta \ml{A}^{\mm{T}}_0 B^0.
\end{eqnarray}
Here and in what follows, the notation $f_0$ (except for the notations $\mathcal{E}_0$, $\mathcal{D}_0$, $\varphi_0$ and $\psi_0$ in \eqref{202111031538}, \eqref{2022111091024} and \eqref{2024}) as well as $f^0$ denotes the value of the function $f$ at $t=0$.
If we assume that the frozen condition holds, i.e.
\begin{equation}
\label{201903081437}
\ml{A}^{\mm{T}}_0 B^0=\bar{M}  \mbox{ (i.e., }B^0= \partial_{\bar{M}}\zeta^0\mbox{)},
 \end{equation}
 where $\bar{M}$ is defined by \eqref{202111031455}, then
 \eqref{0124} reduces to
 \begin{eqnarray}
\label{0124xx}   B=\varpi\partial_\omega\zeta  .
\end{eqnarray}
Here we should point out that $B$ given by \eqref{0124xx} automatically satisfies \eqref{01dsaf16asdfasfdsaasf}$_3$
and \eqref{01dsaf16asdfasfdsaasf}$_5$. Moreover, from \eqref{0124xx} we see that the magnetic tension in Lagrangian coordinates
has the relation
\begin{equation}\nonumber
 B\cdot \nabla_{\ml{A}} B= \varpi^2\partial_{\omega}^2\zeta.
 \end{equation}

Let $I$ denote a $3\times 3$ identity matrix, $m^2= \lambda  \varpi^2 /4\pi\rho   $ and $\eta=\zeta-y$.
Consequently, under the assumption \eqref{201903081437}, the system \eqref{01dsaf16asdfasfdsaasf} is equivalent to
the following system of evolution equations for $(\eta,u,q)$:
\begin{equation}\label{01dsaf16asdfasf}
                              \begin{cases}
\eta_t=u , \\[1mm]
 u_t+\nabla_{\ml{A}}q- \nu \Delta_{\ml{A}}u=
  m^2 \partial_\omega^2\eta, \\[1mm]
\div_\ml{A}u=0,
\end{cases}
\end{equation}
where $B=m\partial_\omega(\eta+y)$
and $ {\mathcal{A}}=(\nabla\eta+I)^{-\mm{T}}  $.

For the well-posedness of \eqref{01dsaf16asdfasf} defined in $\mathbb{T}^3$, we impose the initial condition:
\begin{equation}\label{01dsaf16asdfasfsaf} (\eta, u)|_{t=0}=(\eta^0,  u^0)\mbox{ in } \mathbb{T}^3.
\end{equation}
\subsection{Notations}
 Before stating our main results, we introduce some notations which will be frequently used throughout this paper.
\begin{enumerate}
  \item[(1)] Basic notations: $\langle t\rangle:=t+1 $, $I_T:=(0,T)$ for $0<T\leqslant \infty$, $\overline{I_T}$ is the closure of $I_T$ (in particular, $\overline{I_\infty}=\mathbb{R}^+_0:=[0,\infty)$), $\tilde{\mathcal{A}} :=\mathcal{A}-I $, $\Omega_T:=\mathbb{T}^3\times I_T$,
$\int:=\int_{(-1,1)^3}$, $(w)_{\mathbb{T}^3}:=\int w\mm{d}y$, $\alpha=(\alpha_1,\alpha_2,\alpha_3)$ denotes the multi-index with respect to the variables $y$. In addition,
 \begin{align}
 \sigma(s):=
 \begin{cases}
s &\mbox{for }s\geqslant 0; \\
0 &\mbox{for }s<0 .
 \end{cases} \label{2022111161031}
 \end{align}
  \item[(2)] Simplified Banach spaces:
\begin{align}
& L^r:=L^r (\mathbb{T}^3)=W^{0,r}(\mathbb{T}^3),\;\; {H}^i:=W^{i,2}(\mathbb{T}^3) ,\nonumber \\
&  H^i_\sigma:=\{u\in H^i~|~\mm{div}u=0\},\
  H^{i+1}_{1}:=\{\eta\in H^{i+1} ~|~\det\nabla(\eta +y )=1\},\nonumber \\
 &\mathcal{U}_{T}:=\{u\in C^0(\overline{I_T}, H^{17})\cap L^2({I_T}, H^{18}) ~|~\partial_t^9 u\in L^2(I_T, L^2),
\nonumber \\
&\qquad \qquad \partial_t^i u\in C^0(\overline{I_T}, H^{17-2i})\cap L^2(I_T,H^{18-2i})\mbox{ for }1\leqslant i\leqslant 8\},\nonumber\\
&\underline{\mathcal{U}}_T:=\{u\in \mathcal{U}_{T}~|~(u)_{\mathbb{T}^3}=0\},  \nonumber \\
  &{ \underline{X}:=\{w\in X\cap L^2~|~(w)_{\mathbb{T}^3}=0\}, } \nonumber
  \end{align}
 {where $X$ denotes a Banach space, $1< r\leqslant \infty$ and $i \geqslant 0$ is an integer. }
 \item[(3)] Simplified function classes:
\begin{align}
& H^{18}_{*}:=\{\xi\in H^{18}~|~ \xi(y)+y   : \mathbb{R}^3 \to \mathbb{R}^3 \mbox{ is a } C^{16}\mbox{-diffeomorphism mapping}\},\nonumber\\
&
  \underline{H}^{18,*}_{1,T}:=\{\eta\in C^0(\overline{I_T}, \underline{H}^{18}_{1})~|~ \eta(t) \in H^{18}_{*}\mbox{ for each }t\in \overline{I_T} \}.\nonumber
\end{align}
\item[(4)] Simplified norms and functionals: for integers $i\geqslant 0$,
\begin{align}
&\|\cdot \|_{i} :=\|\cdot \|_{H^i(\mathbb{T}^3)},\ \|\nabla^i \cdot\|^2_0:=\sum_{|\alpha|=i}\| \partial^\alpha \cdot\|^2_0,  \nonumber \\
                 &  \mathcal{E}_{i} := \|(\nabla \eta, u,m\partial_\omega\eta) \|_{i}^2, \  \mathcal{D}_{i} := \|(\nabla u,m\partial_\omega\eta) \|_{i}^2,\label{202111031538} \\
                                                   &  \mathcal{E}_{H}:= \mathcal{E}_{16}
                                              +\|\eta\|_{18}^2+  m^{-2/3} \|( u,m\partial_\omega\eta) \|_{17}^2, \label{2538} \\
   &  \mathcal{D}_H:= \| (u,m\partial_\omega\eta)\|_{17}^2
 + m^{-2/3}\|  u \|_{18}^2   .  \nonumber
    \end{align}
               In addition, we define the following two parameters, which depend  on $m$ and the initial energy functionals:
               \begin{align}
&               \Xi:= \sum_{i=0}^3 (1+m^{-2})^i   \mathcal{E}_{4i}^0 ,\label{20202111110854} \\
 &              \vartheta=( \mathcal{E}_{12}^0)^{1/8} ( 1+(\mathcal{E}_{4}^0 )^{3/2} )\Xi^{3/8}+\mathcal{E}_{6}^0+ (\mathcal{E}_{6}^0)^{2} .
               \label{20211051042}
               \end{align}
\item[(5)]  {General} constants: $c$ and ${c}_i$ ($1\leqslant i \leqslant 3$)
denote constants, which depend on $\omega$ and $\nu$ at most (independent of $m$); moreover
  ${c}$ may vary  from line to line, but ${c}_i$ are fixed.
  $c_0$  denotes a generic constant independent of any parameter.
In addition, we use the notation    $ c({\chi_1,\ldots,\chi_n})$  to denote a generic constant, which only depends on the parameters $\chi_1$, $\ldots$, $\chi_n$.
 Finally, $A\lesssim_0 B$, $A\lesssim B$   and $A\lesssim_{\chi_1,\ldots,\chi_n} B$  mean that $A\leqslant c_0B$, $A\leqslant cB$ and $A\leqslant c({\chi_1,\ldots,\chi_n}) B$, resp..
\end{enumerate}

\subsection{Existence of large perturbation solutions}
Now we state the firs result for the existence result of solutions with large initial perturbation in some class for the initial value  problem \eqref{01dsaf16asdfasf}--\eqref{01dsaf16asdfasfsaf}.
\begin{thm}
\label{201904301948}
If the unit vector $\omega\in \mathbb{R}^3$ satisfy the Diophantine condition \eqref{20201102121755},
then there are positive constants $c_1\geqslant 4$, $c_2>0$ and a sufficiently small constant $c_3\in (0,1]$, such that for any
$(\eta^0,u^0)\in (\underline{H}^{18}_{1}\cap H^{18}_*)\times \underline{H}^{17}$ and $m$ satisfying the incompressible condition
$\mm{div}_{\mathcal{A}^0}u^0=0$  and the condition of strong magnetic field
\begin{equation}
\label{201909281832}m^{-1}\max\left\{\left( c_1 \mathcal{E}_{H}^0 e^{ c_2 \vartheta}\right)^{1/2},
c_1 \mathcal{E}_{H}^0 e^{c_2 \vartheta}\right\}\leqslant c_3 ,
\end{equation}
the initial value problem \eqref{01dsaf16asdfasf}--\eqref{01dsaf16asdfasfsaf} admits a unique global classical solution
$(\eta, u,q)\in \underline{H}^{18,*}_{1,\infty}\times \underline{\mathcal{U}}_\infty\times C^0(\mathbb{R}^+_0,\underline{H}^{17})$.
Moreover, the solution $(\eta,u)$ enjoys, for any $t\geqslant 0$,
\begin{enumerate}[(1)]
  \item the decay-in-time of lower-order derivatives
  \begin{align}
& \sum_{i=0}^{3} \left((1+m^{-2})^i \langle t\rangle^{ (3-i)} {\mathcal{E}}_{4i}+ (1+m^{-2})^i\int_0^t\langle \tau\rangle^{ (3-i)}\mathcal{D}_{4i}  \mm{d}\tau\right)\lesssim \Xi.
  \label{20190safd5041053}
\end{align}
  \item the stability estimates
\begin{align}
& \mathcal{E}_{i} +  \int_0^t \mathcal{D}_{i}\mm{d}\tau\lesssim\mathcal{E}^{0}_{i}\mbox{ for }0\leqslant i\leqslant 12 , \label{2053}\\
& \mathcal{E}_{H} \lesssim\mathcal{E}^{0}_{H}
e^{c_2 \vartheta}  , \label{202211111} \\
&   \int_0^t \mathcal{D}_{H}\mm{d}\tau\lesssim\mathcal{E}^{0}_{H}
(1+\vartheta e^{c_2 \vartheta}).\label{201041053}
\end{align}
\end{enumerate}
In addition, $(\eta,u,q)$ satisfies the additional estimates \eqref{200s12}--\eqref{20200sasf12}, \eqref{2017020614181721} and
\begin{equation}
\label{20202111005120114}
 \sup_{0\leqslant t <\infty} \| \eta\|_{15}  \lesssim_0 1.
\end{equation}
\end{thm}
\begin{rem}
\label{202110121947}
For any given $\tau>N-1$ and $N\geqslant 2$, for almost all $\omega\in \mathbb{R}^N$, there exists a positive constant $c(\omega,N,\tau)$ such that
$$|\chi\cdot\omega|\geqslant c(\omega,N,\tau)|\chi|^{-\tau}\mbox{ for any }\omega\in \mathbb{Z}^N\setminus\{0\}, $$
please refer  to the lemma on p. 139 in \cite{alinhac2007pseudo}.
Using this fact, we easily see that, for any given $\tau>N-1$ and $N\geqslant 2$, for almost all  $\omega\in \partial B$ with respect to the $(N-1)$-dimensional measure,
 where $\partial B:=\{y\in \mathbb{R}^N~|~|y|=1\} $, there exists a constant $c(\omega,N,\tau)$ such that
$$|\chi\cdot \omega|\geqslant c(\omega,N,\tau)|\chi|^{-\tau}\mbox{ for any }\chi\in \mathbb{Z}^N\setminus \{0\}.$$
Hence $\omega$ mentioned in the above Theorem  \ref{201904301948} exists.
\end{rem}
\begin{rem}
The direction condition for $\omega$, such as the Diophantine condition, is necessary in Theorem \ref{201904301948}.
In fact, by the magnetic inhibition mechanism, we do not expect that Theorem \ref{201904301948} can be extended to the case $\omega=e_i$, unless additional structural conditions are imposed. Here $e_i\in \mathbb{R}^3$ is the unit vector with the $i$-th component being $1$.
\end{rem}
\begin{rem}We can easily construct a family of initial data $(\eta^0,u^0)$  {satisfying} all the assumptions in Theorem \ref{201904301948},
where $\eta^0\neq 0$ and $u^0\neq 0$. In fact, let
\begin{align}
\bar{\eta}=\bar{u}=(\sin x_1 \cos x_2 \cos x_3, \cos x_1\sin x_2 \cos x_3, -2\cos x_1\cos x_2 \sin x_3). \nonumber
 \end{align}Because of $\mm{div}\bar{\eta}=\mm{div}\bar{u}=0$, for sufficiently
small $\varepsilon$, there exists a function pair $(\eta^0 ,u^0 )$ enjoying the form
$(\eta^0 ,u^0 )= (\varepsilon\bar{\eta}+\varepsilon^2 \eta^{\mm{r}},\bar{u}+\varepsilon u^{\mm{r}})$,
where $(\eta^{\mm{r}},u^{\mm{r}})\in \underline{H}^{18}\times\underline{H}^{17} $ satisfies
$\| \eta^{\mm{r}}\|_{18}+\| u^{\mm{r}}\|_{17}\leqslant c_0 $,  \begin{align}
   \begin{cases}
   -\Delta \eta^{\mm{r}} +\nabla \beta_1=0, \\
   \mm{div} \eta^{\mm{r}}= \varepsilon^{-2}r_{\varepsilon\bar{\eta}+\varepsilon^2 \eta^{\mm{r}}} , \\
     (\eta^{\mm{r}})_{\mathbb{T}^3}  = 0
   \end{cases}
\end{align}
and
\begin{align} \begin{cases}
-\Delta u^{\mm{r}} +\nabla \beta_2=0, \\
\mm{div} u^{\mm{r}}=-\varepsilon^{-1}\mm{div}_{\tilde{\mathcal{A}}^0}(\bar{u}+\varepsilon u^{\mm{r}}),\\
({u}^{\mm{r}})_{\mathbb{T}^3}  = 0,
\end{cases}
\nonumber
\end{align}
where $r_{\varepsilon\bar{\eta}+\varepsilon^2 \eta^{\mm{r}}}$ is defined as $r_{\eta}$ in \eqref{202108271645} with $\varepsilon\bar{\eta}+\varepsilon^2 \eta^{\mm{r}}$ in place of $\eta$ and $\tilde{\mathcal{A}}^0:=(\nabla \eta^0+I)^{-\mm{T}}-I$.
  We refer the reader to \cite[Proposition 5.1]{JFJSWZhangwei} for a proof. Moreover,
it is easy to check that for sufficiently small $\varepsilon\in (0,1]$, $(\eta^0 ,u^0 )$ constructed above is non-zero,
belongs to $ (\underline{H}^{18}_{1}\cap H^{18}_*)\times \underline{H}^{17} $, and satisfies $\mm{div}_{\mathcal{A}^0}u^0=0$.
We further take $m=\varepsilon^{-1}$ to immediately see that $(\eta^0 ,u^0 )$ and $m$ satisfy \eqref{201909281832}
for sufficiently small $\varepsilon$. Furthermore, $\mathcal{E}_H^0\leqslant c_0$
for some constant $c_0$ independent of $\varepsilon$ and $m$.
\end{rem}
\begin{rem}
\label{20221110152006}
In the above theorem, we have assumed $(\eta^0)_{\mathbb{T}^3}=(u^0)_{\mathbb{T}^3}=0$.
If $((\eta^0)_{\mathbb{T}^3}, (u^0)_{\mathbb{T}^3})\neq 0$, we can define $\bar{\eta}^0:=\eta^0-(\eta^0)_{\mathbb{T}^3}$ and
$\bar{u}^0:=u^0-(u^0)_{\mathbb{T}^3}$. Then, by virtue of Theorem \ref{201904301948}, there exists a unique global classical solution $(\bar{\eta},\bar{u},\bar{q})$ to the initial value problem \eqref{01dsaf16asdfasf}--\eqref{01dsaf16asdfasfsaf} with initial data $(\bar{\eta}^0,\bar{u}^0)$. It is easy to verify that
$(\eta,u,q):=(\bar{\eta}+ t(u)_{\mathbb{T}^3}+ (\eta^0)_{\mathbb{T}^3},\bar{u}+(u)_{\mathbb{T}^3},\bar{q})$
is just the unique classical solution of \eqref{01dsaf16asdfasf}--\eqref{01dsaf16asdfasfsaf} with initial data $({\eta}^0,{u}^0)$.
\end{rem}
\begin{rem}
Recalling the derivation of the decay-in-time \eqref{20190safd5041053}, we easily see that,  higher  the regularity of the initial data is, quicker the decay-in-time of lower-order derivatives of solutions.
Hence we also further establish  a result of almost exponential decay-in-time as in  \cite{TZWYJGw}.
\end{rem}
\begin{rem}
\label{2021110321826}
Let $m\geqslant 1$, $W=B -\bar{M}$, $\Theta_i=\|u \|_{i}^2+\|W\|_{i+4}^2 $ and $\eta^0$ further belongs to $H^{21}$.
It should be noted that the  perturbation magnetic field $W$ is equal to
$\varpi\partial_\omega \eta$ or $2m\sqrt{\pi\rho /\lambda} \partial_\omega \eta$, then we have
$$\mathcal{E}_H^0 \lesssim_{\rho,\omega} \Theta_{17}^0\mbox{ and }\vartheta \leqslant c_4 (\rho,\omega) (\sqrt{\Theta_{12}^0}+ (\Theta_{12}^0)^2) , $$
where the constant $c_4 (\rho,\omega)$ depends on $\rho$ and $\omega$.
Thus we easily  see that, if
 \begin{equation}
\label{2012}
\frac{\max\left\{\sqrt{\mathcal{I}^0 },\mathcal{I}^0 \right\}}{\varpi}\ll  1,
\end{equation}
where $\mathcal{I}^0:= \Theta_{17}^0 e^{ c_2 c_4(\sqrt{\Theta_{12}^0}+ (\Theta_{12}^0)^2)}$,
then the condition of strong magnetic field \eqref{201909281832} can be satisfied. Hence we have mathematically verified Assertion \ref{2020111031231} \emph{in Lagrangian coordinates}.
Moreover  we see from \eqref{2012} that the field intensity $\varpi$ increasingly depends on the $H^{17}\times H^{21}$-norm of the initial perturbation of both the velocity and magnetic field.
\end{rem}
\begin{rem}
Recalling that the solution $\eta$ in Theorem \ref{201904301948} satisfies
\begin{align}
&\zeta:= \eta(y, t)+y   : \mathbb{R}^3\to\mathbb{R}^3 \mbox{ is a }C^{16}\mbox{ diffeomorphism mapping}, \label{2312018031adsadfa21601xx}
\end{align}
we can easily  obtain the existence result of large perturbation solutions  for the original 3D system \eqref{1.1xx} by an inverse transformation of Lagrangian coordinates \cite{JFJSWZhangwei}, and thus provides a  mathematical evidence to support Assertion \ref{2020111031231} for the 3D case. Interesting readers can refer to \cite[Theorem 2.3]{JFJSAAD} for writing down the existence result of large perturbation solutions  of the 3D system \eqref{1.1xx}.
\end{rem}

Now we briefly describe the proof idea of Theorem \ref{201904301948}.
We rewrite \eqref{01dsaf16asdfasf}$_2$--\eqref{01dsaf16asdfasf}$_3$ as a nonhomogeneous system:
\begin{equation}\label{s0106pnnnn}
 \begin{cases}
u_t-\nu\Delta  u - m^2\partial_\omega^2\eta =\mathfrak{N} , \\[1mm]
\div u ={{-}}\mm{div}_{\tilde{\mathcal{A}}}u,
\end{cases}
\end{equation}
where $ \mathfrak{N}:=\mathcal{N}^\nu-\nabla_{\tilde{\ml{A}}}q-\nabla q$,
$\mathcal{N}^\nu:= \partial_l(\mathcal{N}^\nu_{1,l}, \mathcal{N}^\nu_{2,l}, \mathcal{N}^\nu_{3,l})^{\mm{T}}$ and
$\mathcal{N}^\nu_{j,l} := \nu (\mathcal{A}_{kl}\tilde{\mathcal{A}}_{km} +\tilde{\mathcal{A}}_{ml} )\partial_mu_j$. It should be noted that the linearized pressure term $\nabla q$ can be regarded as a nonlinear term under \emph{the sense of energy integrals}.

If we ignore the nonlinear terms in \eqref{s0106pnnnn}, then by a method of \emph{a priori} estimates,
\begin{align}
\|(\nabla \eta, u, m\partial_{\omega}\eta)\|_{i+4}\lesssim  \|(\nabla \eta^0,u^0, m\partial_{\omega}\eta^0)\|_{i+4}.
\label{20202111032017}
\end{align}
 Thus, exploiting the Diophantine condition \eqref{20201102121755} and the generalized  Poinc\'are's inequality \eqref{2020211013200910}, we  obtain
\begin{align}
\label{202211132010}
\|\nabla \eta\|_i\lesssim  \|(\nabla \eta^0, u^0, m\partial_{\omega}\eta^0)\|_{i+4}/m.
\end{align}
Motivated by the approximate theory of non-resistive MHD equations in \cite{BCSCSPLL,ZTGS}, we naturally guess that \eqref{202211132010} shall hold in the nonlinear equations \eqref{s0106pnnnn} for sufficiently large $m$.

 In turn, if  \eqref{202211132010} exists in the nonlinear equations \eqref{s0106pnnnn}, it is easy to see the system \eqref{s0106pnnnn}  can be approximated by the corresponding linearized system. Since the linear system admits a global solution, the nonlinear system \eqref{s0106pnnnn} may also admit a global solution
with  large data under a strong magnetic field.
Thus, motivated by \eqref{20202111032017}, for given initial energy  $\mathcal{E}_{17}^0$, we naturally expect  to derive the \emph{a priori} estimate of $(\eta,u)$ like
\begin{align}\nonumber
\|(\nabla \eta, u, m\partial_{\omega}\eta)\|_{17}\leqslant  K/2 \mbox{ for sufficiently large } m
\end{align}
under the \emph{a priori} assumption
 \begin{align}
\nonumber
\|(\nabla \eta, u, m\partial_{\omega}\eta)\|_{17}\leqslant K\mbox{ with }Km^{-1}\ll 1.
\end{align}

However, the nonlinear terms appearing in \eqref{s0106pnnnn} destroy the above expectation. In fact, when we perform the estimates of highest-order derivatives, an   integral term increasing on $m$ appears, see the last term in \eqref{2022111040851}. To balance the increasing term, we shall adjust the highest-order energy functional. This is the reason why the structures of the both energy and dissipation of highest-order in \eqref{2538} are different to the ones of lower-order in \eqref{202111031538}. After this adjustment, we can use a three-energy method to conclude that there are constants $K$ (increasingly depending on $\vartheta$ and $ \mathcal{E}_{H}^0$) and $\delta$, such that
\begin{equation}
\label{aprpiose1}
\sup_{0\leqslant t\leqslant T} ( \|\eta(t)\|_{18}^2+ m^{-2/3} \|  u  \|_{17}^2+
 \|(u,m\partial_\omega \eta)(t)\|_{16}^2
)
   \leqslant{K}^2/ 4,
\end{equation}
if
\begin{align}
&\sup_{0\leqslant t\leqslant T}( \|\eta(t)\|_{18}^2+
 \|(u,m\partial_\omega \eta)(t)\|_{16}^2
)
\leqslant K^2 \mbox{ for any given } T>0\label{aprpiosesnewxxxx}
\end{align}
and
\begin{equation}
\label{aprpiose1snewxxxxz}
\max\{K  , K^2\}/m\in (0,\delta]\mbox{ with }\delta\ll 1.\end{equation}

The  \emph{a priori} stability estimate \eqref{aprpiose1}, together with a local well-posedness result on \eqref{01dsaf16asdfasf}--\eqref{01dsaf16asdfasfsaf},
immediately yields Theorem \ref{201904301948}.  Here we explain how to perform the three-energy method, which had been widely used in the investigation of the problems involving wave phenomena, see \cite{GYTIAE2,GYTIDAP} for examples. Roughly speaking, we call \eqref{20190safd5041053}, \eqref{2053} with $i=12$ and \eqref{202211111} the lower-order,  higher-order and highest-order energy inequalities, resp. Under the assumptions \eqref{aprpiosesnewxxxx} and
\eqref{aprpiose1snewxxxxz}, we first get the higher-order energy inequality, and then further obtain the lower-order energy inequality. Thanks to the decay-in-time of lower-order energy, finally  we can close the highest-order energy inequality, see Section \ref{20190sdafs4301948} for the detailed performance.

\subsection{Vanishing phenomena of nonlinear interactions}
Now we turn to mathematically stating the vanishing phenomena of nonlinear interactions with respect to $t$ and $m$.
 \begin{thm}\label{201912041028}
Let the global solution of \eqref{01dsaf16asdfasf}--\eqref{01dsaf16asdfasfsaf} $(\eta,u,q)$ be given by Theorem \ref{201904301948}.
\begin{enumerate}[(1)]
   \item Then, we can use the initial data of $(\eta,u)$ to construct a function pair
   $(\eta^{\mm{r}}, u^{\mm{r}})\in  \underline{H}^{18}\times \underline{H}^{17}$, such that the following linear pressureless
   initial value problem
   \begin{equation}\label{202001070914}
 \begin{cases}
\eta_t^{\mm{L}}=u^{\mm{L}} , \\[1mm]
u_t^{\mm{L}} -\nu\Delta  u^{\mm{L}}= m^2\partial_\omega^2\eta^{\mm{L}}  , \\[1mm]
\div u^{\mm{L}} =0,\\
(\eta^{\mm{L}},u^{\mm{L}})|_{t=0}=(\eta^0+\eta^{\mm{r}}, u^0+u^{\mm{r}}) \in
\underline{H}^{18}_\sigma\times \underline{H}^{17}_\sigma
\end{cases}
\end{equation}
admits a unique classical solution
$(\eta^{\mm{L}}, u^{\mm{L}} )\in C^0(\mathbb{R}_0^+, \underline{H}^{18}_\sigma)\times\underline{\mathcal{U}}_\infty $. Moreover, the linear solution
$(\eta^{\mm{L}}, u^{\mm{L}} )$ enjoys the following estimates
\begin{align}
& \sum_{i=0}^{3} \left((1+m^{-2})^i \langle t\rangle^{ (3-i)} {\mathcal{E}}_{4i}^{\mm{L}}+ (1+m^{-2})^i\int_0^t\langle t\rangle^{ (3-i)}\mathcal{D}_{4i}^{\mm{L}}  \mm{d}\tau\right)\lesssim \sum_{i=0}^3 (1+m^{-2})^i   {\mathcal{E}}_{4i}^{\mm{L}}|_{t=0}   , \label{2020052010sdfa33}\\
 &  \mathcal{E}_{j}^{\mm{L}} +
\int_0^t\mathcal{D}_{j}^{\mm{L}} \mm{d}\tau\lesssim \mathcal{E}^{\mm{L} }_{j}|_{t=0} , \label{202033}
\end{align}
where $0\leqslant j\leqslant 17$ and $$\mathcal{E}^{\mm{L}}_\chi
 := \|(\nabla \eta, u,m\partial_\omega\eta) \|_{\chi}^2 \mbox{ and  }  \mathcal{D}^{\mm{L}}_\chi := \|(\nabla u, m\partial_\omega\eta) \|_{\chi}^2 . $$
In addition  the function pair $(\eta^{\mm{r}}, u^{\mm{r}}) $ satisfies
\begin{align}
&\mm{div}(u^0+u^{\mm{r}})=\mm{div}(\eta^0+\eta^{\mm{r}})=0,\label{202011032122}  \\
  &\| u^{\mm{r}} \|_{k}\lesssim_0   \| \eta^0\|_3\| u^0\|_{k}+
  \|\eta^0\|_{k}\| u^0\|_3,\label{202af011032122}
  \\
  & \|\partial_\omega \eta^{\mm{r}}\|_k\lesssim_0 \|\eta^0\|_3 \|\partial_\omega \eta^0\|_k+ \|\partial_\omega  \eta^0\|_3\|\eta^0\|_k,\ \|\eta^{\mm{r}}\|_{l}\lesssim_0  \|\eta^0\|_3\|\eta^0\|_{l}, \label{202011032123}
\end{align}
where  $1\leqslant k\leqslant 17$ and   $1\leqslant l\leqslant 18$.
 \item  Let $(\eta^{\mm{d}}, u^{\mm{d}} )= (\eta-\eta^{\mm{L}} , u-u^{\mm{L}} ) $, then for any $t\geqslant 0$,
 \begin{align}
& \mathcal{E}_{12} ^{\mm{d}}  +  \int_0^t \mathcal{D}_{12} ^{\mm{d}} \mm{d}\tau
 \lesssim m^{-1}  \mathcal{E}^{0}_{17}
e^{c_2 \vartheta} \left( \sqrt{\Xi+\mathcal{E}^{0}_{12}} + \mathcal{E}^{0}_{12} \right), \label{2020sfa05201033} \\
& \sum_{i=0}^{3} \left((1+m^{-2})^i \langle t\rangle^{ (3-i)} {\mathcal{E}}_{4i}^{\mm{d}}+ (1+m^{-2})^i\int_0^t\langle \tau\rangle^{ (3-i)}\mathcal{D}_{4i}^{\mm{d}}\mm{d}\tau  \right)  \lesssim m^{-1}\left( {\mathcal{E}_{15}^0}\Xi+ \Phi\right)
\label{201905041053xxx},
\end{align}
\end{enumerate}
where
\begin{align}
& \Phi:=
 \Xi  \left(\sqrt{\mathcal{E}^{0}_{12}} +  \mathcal{E}^{0}_{12}
 \right)\sum_{i=0}^{2} (1+m^{-2})^i+  \mathcal{E}^{0}_{17}
e^{c_2 \vartheta} \left( \sqrt{\Xi+\mathcal{E}^{0}_{12}} + \mathcal{E}^{0}_{12} \right)(1+m^{-2})^3, \label{2020211111342}\\
&\mathcal{E}^{\mm{d}}_\chi:=\|(\nabla \eta^{\mm{d}}, u^{\mm{d}},m\partial_\omega\eta^{\mm{d}}) \|_{\chi}^2\mbox{ and }\mathcal{D}_{\chi}^{\mm{d}} := \|( \nabla u^{\mm{d}},m\partial_\omega\eta^{\mm{d}}) \|_{\chi}^2. \nonumber
\end{align}
\end{thm}
\begin{rem}
Exploiting the interpolation inequality \eqref{202108261405}, we easily derive  from \eqref{2020sfa05201033} and \eqref{201905041053xxx} that, for any $0\leqslant i\leqslant 12$,
$$
\mathcal{E}^{\mm{d}}_i \lesssim   m^{-1}\left(\mathcal{E}_{15}^0\Xi+ \Phi \right)^{(12-i)/12} \left( \mathcal{E}^{0}_{17}
e^{c_2 \vartheta} \left( \sqrt{\Xi+\mathcal{E}^{0}_{12}} + \mathcal{E}^{0}_{12} \right) \right)^{i/12}.
$$
In addition, for $m\geqslant 1$, the estimate \eqref{201905041053xxx} can be simplified  as follows:
$$\sum_{i=0}^{3} \left( \langle t\rangle^{ (3-i)} {\mathcal{E}}_{4i}^{\mm{d}}+ \int_0^t\langle \tau\rangle^{ (3-i)}\mathcal{D}_{4i}^{\mm{d}} \mm{d}\tau \right)  \lesssim   m^{-1}{\mathcal{E}^{0}_{17}
e^{c_2 \vartheta}}\left( \sqrt{\mathcal{E}_{12}^0}+ \mathcal{E}_{12}^0\right) .$$
\end{rem}
\begin{rem}
The decay-in-time in \eqref{2020052010sdfa33} can be easily observed, since the both linear and nonlinear solutions enjoy the same decay-in-time, see \eqref{20190safd5041053} and \eqref{2020052010sdfa33}.
Noting that the inhomogeneous term $\mathfrak{N}$
 includes $\partial_i\eta_j$ (the term $\nabla q$ should be understood under the sense of the energy integrals) with $1\leqslant i$, $j\leqslant 3$, and
\begin{align}
\label{20202111041459}
\|\nabla \eta\|_{L^\infty}\lesssim_0 \|\partial_\omega \eta\|_6\leqslant K/2m
\end{align}
by \eqref{aprpiose1}, \eqref{2020211013200910} and \eqref{202004221saffad412}, we formally see from the error equations in \eqref{01dsaf16safafasdfasfxx} the appearance of $m^{-1}$ in \eqref{2020sfa05201033} and \eqref{201905041053xxx}.
The two estimates \eqref{2020sfa05201033} and \eqref{201905041053xxx} present    the physical phenomena of vanishing of nonlinear interactions  for large time and strong magnetic field.
 \end{rem}

  We can follow the idea of deriving the estimates \eqref{20190safd5041053} and \eqref{2053} to establish Theorem \ref{201912041028}, the proof of which will be presented in Section \ref{2020011sdf92326}.
Here we explain why we have to modify $(\eta^0,u^0)$
to be initial data of $(\eta^{\mm{L}},u^{\mm{L}} )$ in \eqref{202001070914}$_4$.
\begin{enumerate}
  \item[(1)]
Since the initial data $u^{\mm{L}}|_{t=0}$ is divergence-free, i.e.  $\mm{div}(u^{\mm{L}}|_{t=0})=0$, we have  to adjust
the initial data $u^0$ as in \eqref{202001070914}$_4$.
  \item[(2)]
 { If the initial data $\eta^0$ of $\eta$ is directly used to be initial data of the corresponding linear problem, then
 we see that $\mm{div}\eta^{\mm{L}}=\mm{div}\eta^0$, and a time-decay of $\partial_\omega\eta^{\mm{d}}$ in \eqref{201905041053xxx}
 can not be expected unless $\mm{div}\eta^0=0$. Hence, we have to modify $\eta^0$ as in \eqref{202001070914}$_4$,
  so that the obtained new initial data ``$\eta^0+\eta^{\mm{r}}$'' is also divergence-free. }
\end{enumerate}

Finally we mention that recently  some authors studied the case of inviscid, non-resistive MHD fluids with velocity damping,
i.e., the viscosity term in the system \eqref{1.1xx} is replaced by the velocity damping term $\kappa \rho v$ with $\kappa$ being
the damping coefficient  \cite{WJHWYFXXJG,DYYWZYOSJMA,JFJSAAD}.
Following the arguments of Theorems \ref{201904301948}, \ref{201912041028} and \cite[Theorem 2.5]{JFJSAAD}, we can extend our aforementioned results in  Theorems \ref{201904301948} and \ref{201912041028} to the inviscid case with the damping term $\kappa\rho v$, and show that for the inviscid
case with damping, the decay in time  is exponentially fast;
while the convergence rate in $m$ as $m\to\infty$ of a classical solution of the original nonlinear system to the solution
of the corresponding linear system is in the form of $m^{-1}$, which is faster than $m^{-1/2}$ in the viscous case.

\section{Proof of Theorem \ref{201904301948} } \label{20190sdafs4301948}

This section is devoted to the proof of Theorem \ref{201904301948}. First we derive some energy (\emph{a priori}) estimates for the solution $(\eta,u)$ of the    initial value problem \eqref{01dsaf16asdfasf}--\eqref{01dsaf16asdfasfsaf}  under the \emph{a priori} assumption \eqref{aprpiosesnewxxxx} associated with the relative smallness condition \eqref{aprpiose1snewxxxxz} in Subsection \ref{202061410014}, then further establish the desired \emph{a priori} estimate \eqref{aprpiose1} in Subsection \ref{202006141016}, and finally, introduce a local well-posedness result for \eqref{01dsaf16asdfasf}--\eqref{01dsaf16asdfasfsaf} and complete the proof of
 Theorem \ref{201904301948} by a standard continuity argument in Subsection \ref{20206141015}. In the derivation of \eqref{aprpiose1}, we obtain the stability estimates \eqref{2053}--\eqref{201041053} as by-products.

\subsection{Energy estimates} \label{202061410014}

Let $(\eta,u,q)$ be a solution of the initial value problem \eqref{01dsaf16asdfasf}--\eqref{01dsaf16asdfasfsaf} defined on $\Omega_T$  for any given $T>0$, and the initial data $(\eta^0,u^0) $ is a non-zero function, belongs to $\underline{H}^{18}_{1} \times \underline{H}^{17}$, and satisfies  $\mm{div}_{\mathcal{A}^0}u^0=0$. Obviously $(\eta,u)$ automatically satisfies  $(\eta)_{\mathbb{T}^3}=(u)_{\mathbb{T}^3}=0$ as well as the initial data $(\eta^0,u^0)$.
We further assume that $(\eta,u,q)$ and $K$ satisfy  $(q)_{\mathbb{T}^3}=0$ and \eqref{aprpiosesnewxxxx}, where $K $ will be defined by \eqref{201911262060}.

Before deriving the energy estimates for $(\eta,u)$ in Lemmas \ref{qwe201612132242nn}--\ref{qwe201612132242nxsfssdfsxx}, we shall introduce some basic inequalities and establish some preliminary estimates for $(\eta,u)$ by the following three lemmas.

\begin{lem}\label{202005091306}
 We have the following basic inequalities:
\begin{enumerate}[(1)]
\item For any $s_1$, $s_2$, $s_3\in \mathbb{R}$ satisfying $s_1\leqslant s_2\leqslant s_3$,  and any positive constant $a$, it holds that
\begin{align}
 a^{s_2}\leqslant  a^{s_1}+a^{s_3}.
 \label{202021100321520}
 \end{align}
 The above inequality is obvious by the monotonicity of $a^{s}$ with respect to the variables $s$.
   \item
Poinc\'are's inequality \cite{NASII04}: for any given $i\geqslant 0$,
\begin{align}
&\label{20160614fdsa19asfda57x}
 \| f\|_i \lesssim_i \|\nabla^i f\|_0 \mbox{ for any }f\in \underline{H}^i.
\end{align}
   \item
 Generalized  Poinc\'are's inequalities \cite{chen2021}: if $\omega \in \mathbb{R}^3$ satisfies the Diophantine condition \eqref{20201102121755}, then it holds that, for any given $i\geqslant 0$,
\begin{align}
\label{2020211013200910}
\|f\|_i\lesssim_{\omega,i} \|\partial_\omega  f\|_{i+3}\mbox{ for any }f \in \underline{H}^{ 3+i}(\mathbb{T}^{3}).
\end{align}
 \item Interpolation inequalities \cite{ARAJJFF}:   for  given $0\leqslant j \leqslant i$,
\begin{align}
&\|f\|_{L^\infty } \lesssim_0  \|f\|_1^{1/2}\| f \|_2^{1/2}\mbox{  for any }f\in H^2,
\label{202004221saffad412} \\
&\label{202108261405}
\|f\|_{j}\lesssim_{i,j}\|f\|_{0}^{1-j/i}\|f\|_{i}^{j/i}
  \mbox{  for any }f\in H^i.
\end{align}
\item
 Product estimates  \cite{BJAFP315200}:
% \begin{align}d
% &\| fg\|_0lesssim_0  \|f\|_{L^6} \| g\|_{L^3}\lesssim_0  \|f\|_{1} \sqrt{\| g\|_{0}\| g\|_{1}} \mbox{ for any }f,\ g\in   {H}^1, \label{202005safsda08}\\
 \begin{align}
 & \| fg\|_i\lesssim_0
 \|f\|_1\sqrt{\|g\|_{i}\|g\|_{1+i}}  \mbox{ for any }(f,g)\in  H^1\times {H}^{i+1},\label{2020saf05sda08}
  \\
&\| fg\|_j\lesssim_j \|f\|_{L^\infty} \| g\|_j+\|f\|_{j}\|g\|_{L^\infty} \mbox{ for any }f,\ g\in   {H}^j\cap L^\infty, \label{2020sadf05sda08}
%\lesssim_i &\|f\|_{0}^{1/4}\|f\|_{0}^{3/4} \| g\|_i+\|f\|_{i}\|g\|_{0}^{1/4}\|g\|_{2}^{3/4} \mbox{ for any }f,\ g\in   {H}^i\cap L^\infty. \label{202005sda08}
\end{align}
where $0\leqslant i\leqslant 1$ and $j\geqslant 0$ are given.
In particular, the two estimates \eqref{202004221saffad412} and \eqref{2020sadf05sda08} imply that, for any given $i\geqslant 0$,
\begin{align}
\| fg\|_i\lesssim_i  \|f\|_{2} \| g\|_i+\|f\|_{i} \|g\|_{2}  \mbox{ for any }f,\ g\in   {H}^i\cap  H^2. \label{202005sda08}
\end{align}
%\item There is a constant $\delta\in (0,1]$, such that for any $(\xi-y)\in H$ satisfying $\|\nabla (\xi-y)\|_{L^\infty}\leqslant \delta$,\begin{align}\|\nabla_{\mathcal{B}}f\|_0\lesssim_0  \| \nabla f\|_0\lesssim_0 \|\nabla_{{\mathcal{B}}} f\|_0\;\;\mbox{ for any }\nabla f\in L^2,\label{202010261646}\end{align}where $\ml{B}:= \nabla {{\xi}} ^{-\mm{T}}$, see \cite{JFJSJMFMOSERT}.
\end{enumerate}
\end{lem}
\begin{lem}
\label{20202110261955}
Let
$\eta$   further satisfy
 \begin{equation}
\label{aprpi1x}
 \sup_{0\leqslant t \leqslant T} \| \eta\|_3   \lesssim_0    \delta \in (0,1],
\end{equation}  then
\begin{align}
& \|\tilde{\mathcal{A}}\|_{i}  \lesssim_0 \|  \eta \|_{i+1}  \mbox{ for } 0\leqslant i\leqslant 16,\label{200s12}\\
 &  \|\mathcal{A}_t\|_{i}^2 \lesssim_0 \|  u\|_{i+1} +
 \begin{cases}0  &\mbox{for }0\leqslant i\leqslant 2;
\\
  \| \eta\|_{i+1}
\| u\|_3
  & \mbox{for } 3 \leqslant  i\leqslant 16
  \end{cases}\label{20200s12}
\end{align}
and
\begin{align}
 &  \|(\partial_\omega^2 \mm{div} \eta ,\mm{div}_{\tilde{\mathcal{A}}}\partial_\omega^2 \eta )\|_{i}^2 \lesssim_0 F_i
,\label{20200sasf12}
\end{align}
where \begin{align} F_i:=\begin{cases}
\| \partial_\omega   \eta \|_2 \|  \partial_\omega \eta \|_3
    + \| \eta \|_3 \| \partial_\omega^2 \eta \|_2 & \mbox{for }i=1;
  \\
\|  \eta \|_3 \| \partial_\omega^2 \eta \|_{i+1}+ \| \partial_\omega  \eta \|_3 \|  \partial_\omega \eta \|_{i+1}+ \| \eta \|_{i+1}
  (\| \partial_\omega^2 \eta \|_3+   \| \partial_\omega  \eta \|_3 ^2)
    & \mbox{for }10\leqslant i\leqslant 15.
  \end{cases}\nonumber
\end{align}
\end{lem}
\begin{pf}
Noting that $\det(\nabla \eta^0 +I)=1$, then
\begin{align}
\det(\nabla \eta +I)=1. \label{2022110115121}
\end{align}
Thus it is easy to compute out that
\begin{align}
\tilde{\mathcal{A}} = \tilde{\mathcal{A}}^L+ \tilde{\mathcal{A}}^N ,
\label{20201010261550}
\end{align} where
 $$  \tilde{\mathcal{A}}^{{L}} :=\left(\begin{array}{ccc}
\partial_2\eta_2+\partial_3\eta_3 &
-\partial_1\eta_2 &
- \partial_1\eta_3 \\
-\partial_2\eta_1 &
\partial_1\eta_1+\partial_3\eta_3 &
- \partial_2\eta_3 \\
-\partial_3\eta_1 &
- \partial_3\eta_2 &
\partial_1\eta_1+\partial_2\eta_2
        \end{array}\right)
$$
and
$$ \tilde{\mathcal{A}}^{{N}}:=\left(\begin{array}{ccc}
\partial_2\eta_2\partial_3\eta_3-\partial_2\eta_3\partial_3\eta_2 &
\partial_1\eta_3\partial_3\eta_2-\partial_1\eta_2\partial_3\eta_3 &
\partial_1\eta_2\partial_2\eta_3-\partial_1\eta_3\partial_2\eta_2    \\
\partial_2\eta_3\partial_3\eta_1-\partial_2\eta_1\partial_3\eta_3 &
\partial_1\eta_1\partial_3\eta_3-\partial_1\eta_3\partial_3\eta_1 &
\partial_1\eta_3\partial_2\eta_1-\partial_1\eta_1\partial_2\eta_3    \\
\partial_2\eta_1\partial_3\eta_2-\partial_2\eta_2\partial_3\eta_1 &
\partial_1\eta_2\partial_3\eta_1-\partial_1\eta_1\partial_3\eta_2 &
\partial_1\eta_1\partial_2\eta_2-\partial_1\eta_2\partial_2\eta_1
        \end{array}\right).
$$

Applying $\|\cdot\|_i$ and $\|\partial_t\cdot\|_i$ to the identity
\eqref{20201010261550}, resp., and then using \eqref{01dsaf16asdfasf}$_1$, \eqref{202005sda08} and \eqref{aprpi1x}, we immediately get \eqref{200s12} and \eqref{20200s12}.

By \eqref{2022110115121}, we get
\begin{align}\label{202108271645}
\mm{div}\eta= r_{\eta},
\end{align}
where $r_{\eta}:=r_2^{\eta}+r_3^{\eta}$ with
$$\begin{aligned}
&r_2^{\eta}:=\partial_2\eta_1\partial_1\eta_2+
\partial_2\eta_3\partial_3\eta_2+
\partial_3\eta_1\partial_1\eta_3
- \partial_1\eta_1\partial_2\eta_2-\partial_1\eta_1\partial_3\eta_3-
\partial_2\eta_2\partial_3\eta_3
,\\ &r_3^{\eta}:=\partial_1\eta_1(\partial_2\eta_3\partial_3\eta_2
-\partial_2\eta_2\partial_3\eta_3)
+\partial_2\eta_1(\partial_1\eta_2\partial_3\eta_3-\partial_1\eta_3\partial_3\eta_2)
+\partial_3\eta_1(\partial_1\eta_3\partial_2\eta_2-\partial_1\eta_2\partial_2\eta_3).
\end{aligned}$$

Thanks to \eqref{201909261909}, \eqref{20201010261550} and \eqref{202108271645}, we have
\begin{align}
\partial_\omega^2 \mm{div} \eta=  \partial_\omega^2r_2^{\eta}+ \partial_\omega^2r_3^{\eta} \nonumber
\end{align}
and
\begin{align}
\mm{div}_{\tilde{\mathcal{A}}}\partial_\omega^2\eta=
 \tilde{\mathcal{A}}^L  :\partial_\omega^2\nabla \eta +
  \tilde{\mathcal{A}}^N  :\partial_\omega^2 \nabla \eta .\nonumber
\end{align}

Applying $\|\cdot\|_i$ to the two identities above,   and then using \eqref{2020saf05sda08}, \eqref{202005sda08} and \eqref{aprpi1x}, we immediately get \eqref{20200sasf12}.
\hfill $\Box$
\end{pf}

\begin{lem}\label{201612132242nx}
Under the condition \eqref{aprpi1x} with sufficiently small $\delta$, we have
\begin{align}
&\| q\|_{i+2} \lesssim_0
\begin{cases}
\|   u\|_2 \|  u\|_3+ m^2 F_1& \hbox{for }i= 1 ; \\
%    \|   u\|_2 \| u\|_4+ \|u\|_3^2  + m^2F_2&  \hbox{for }i=2 ; \\
 \|   u\|_2 (\| u\|_{i+2}+ \|\eta\|_{i+2}\|u\|_3)  &\\
 + \|  u\|_3 \|   u\|_{i+1}+ m^2 (F_{i} +\| \eta\|_{i+2} F_1)
  &\hbox{for }10\leqslant i\leqslant 15.
 \end{cases} \label{2017020614181721}
\end{align}
\end{lem}
\begin{pf}
(1) By \eqref{201909261909} and \eqref{01dsaf16asdfasf}$_3$, we have
\begin{align} \mm{div}_{\mathcal{A}} u_t = -\mm{div}_{\mathcal{A}_t} u = -\mm{div}({\mathcal{A}}_t^{\mm{T}} u)   \label{202005101008}
\end{align}
and
\begin{align}  \mm{div}_{\mathcal{A}} \nabla_{\mathcal{A}} q  = \mm{div} ({\mathcal{A}}^{\mm{T}}\nabla_{\mathcal{A}}q  )= \mm{div} ( \nabla q+
(\tilde{\mathcal{A}}^{\mm{T}}+\tilde{\mathcal{A}} )\nabla q
+ \tilde{\mathcal{A}}^{\mm{T}}\nabla_{\tilde{\mathcal{A}}}q
). \label{202005101asdfa008}
\end{align}

Let $i=1$, $10$, $\ldots$, $15$ and $\alpha$ satisfy $|\alpha|=i$.
 Applying $\partial^\alpha \mm{div}_{\mathcal{A}}$ to \eqref{01dsaf16asdfasf}$_2$, and then  using \eqref{202005101008}, \eqref{202005101asdfa008} and the fact ``$\mm{div}_{\mathcal{A}}\Delta_{\mathcal{A}}u=0$", we obtain
\begin{equation}
 \Delta  \partial^\alpha  q =   \partial^\alpha f,  \label{202006111729}
\end{equation}
where
$$
\begin{aligned}
f:= &  \mm{div}({\mathcal{A}}_t^{\mm{T}} u)     - \mm{div}( ({\tilde{\ml{A}}}+ {\tilde{\ml{A}}}^{\mm{T}})\nabla q
+{\tilde{\ml{A}}}^{\mm{T}}\nabla_{\tilde{\ml{A}}}q)  + m^2( \partial_\omega^2\mm{div} \eta +\mm{div}_{\tilde{\mathcal{A}}} \partial_\omega^2\eta ) .
\end{aligned}$$

Multiplying \eqref{202006111729} by $\Delta  \partial^\alpha  q$ in $L^2$, and then using regularity theory of elliptic equations, we get
\begin{align}
  \|    \partial^\alpha  q\|_2^2\lesssim_0 \|\Delta  \partial^\alpha  q\|_0^2 \lesssim_0\int|  \partial^\alpha  f  \partial^\alpha  \Delta q|\mm{d} y. \label{201711118xx41x}
\end{align}
The integral term on the right hand side can be  {bounded} as follows:
\begin{align}
\int|   \partial^\alpha  f  \Delta \partial^\alpha  q|\mm{d} y\lesssim_0 &  (\| {\mathcal{A}}_t \|_{i+1}\| u\|_2 +
\| {\mathcal{A}}_t \|_{2}\|  u\|_{i+1}+
 m^2 \| ( \partial_\omega^2\mm{div} \eta +\mm{div}_{\tilde{\mathcal{A}}} \partial_\omega^2\eta ) \|_i  \nonumber \\
& + (1+\|\tilde{\mathcal{A}} \|_{ 2})(\|\tilde{\mathcal{A}}  \|_{i+1}
   \| q\|_3 + \| \tilde{\mathcal{A}} \|_{2}
   \| q\|_{i+2}) ) \| q \|_{i+2}\nonumber \\  \lesssim_0 &
( \|  u\|_2 \| \mathcal{A}_t\|_{i+1} + \| u\|_3 \| u\|_{i+1} + m^2  \| ( \partial_\omega^2 \mm{div} \eta +\mm{div}_{\tilde{\mathcal{A}}} \partial_\omega^2 \eta ) \|_i \nonumber \\
& +  \|  \eta  \|_{i+2}
   \| q\|_3 +\delta
   \|q\|_{i+2} )  \|  q \|_{i+2}, \label{201909121931}
\end{align}
where we have used H\"older's inequality, \eqref{202005sda08} in the first inequality,  and \eqref{aprpi1x}--\eqref{20200s12} in the last inequality,

Making use of the fact $(q)_{\mathbb{T}^3}=0$, Poinc\'are's inequality \eqref{20160614fdsa19asfda57x},  \eqref{aprpi1x}, \eqref{20200s12}  and  \eqref{20200sasf12} we  get \eqref{2017020614181721} from \eqref{201711118xx41x} and \eqref{201909121931} for sufficiently small $\delta$.
  \hfill$\Box$
\end{pf}

From now on, we further assume that $(\eta,u)$ satisfies \eqref{aprpiosesnewxxxx} with \eqref{aprpiose1snewxxxxz}, in which $\delta\in (0,1]$ is a sufficiently small constant. It should be noted that the smallness of $\delta$ only depends on the parameter $\nu$ and the known unit vector $\omega$.
Exploiting  \eqref{2020211013200910} and the interpolation inequality \eqref{202108261405}, we can derive from
\eqref{aprpiosesnewxxxx} and \eqref{aprpiose1snewxxxxz} that
\begin{align}
\label{aprpiosesnewxxxx1x}
& \sup_{0\leqslant t \leqslant T} \| \eta\|_{15}^2  \lesssim_0   \sup_{0\leqslant t \leqslant T}(\| \partial_\omega\eta\|_{16}\|\eta\|_{17}) \lesssim_0 K^2m^{-1}\lesssim_0 \delta \in (0,1],\\
& \sup_{0\leqslant t \leqslant T}
  (\| \eta\|_{13}+ \|\partial_\omega \eta\|_{16} ) (\|u
  %,m\partial_\omega \eta)
\|_{16}+ \|  \eta\|_{18})   \lesssim_0 K^2m^{-1}\lesssim_0 \delta
\label{aprpxx1x}
\end{align}
and
\begin{align}
&  \sup_{0\leqslant t \leqslant T}(
\|\eta\|_{13}
%\| \partial_\omega \eta\|_{16}
       +(\|\eta\|_{18}+\|u\|_{16}^{(2-a)/2}\|m \partial_\omega\eta\|_{16}^{a/2})/m) \lesssim_0 K m^{-1}\lesssim_0 \delta \label{apxxx1x}
\end{align}  for any $0\leqslant a\leqslant 2 $.

Now we proceed to derive some energy estimates for $(\eta,u)$.
\begin{lem}\label{qwe201612132242nn}
Under the conditions \eqref{aprpiosesnewxxxx}--\eqref{aprpiose1snewxxxxz} with sufficiently small $\delta$, we have
\begin{align}
&\frac{\mm{d}}{\mm{d}t} \|\nabla^{i} (u, m\partial_\omega \eta ) \|_0^2   + \nu\| u \|_{i+1}^2
\nonumber \\
&\lesssim \begin{cases}
\delta m^2\|\partial_\omega \eta\|_i^2 &\mbox{for }0\leqslant i \leqslant 12;\\
\delta   m^2  \| \partial_\omega \eta \|_{16}^2 +  \|\eta\|_{18}( \|\eta\|_{17}  (\| u\|_3^2+ \| u\|_3^4& \\
+ m^4    \| \partial_\omega \eta \|_4^2   \|  \partial_\omega   \eta \|_6^2) +m^2( \| \partial_\omega  \eta \|_3 \|\partial_\omega \eta \|_6 \|\eta\|_{16}\|   u\|_3 &\\
+ \|   u\|_3(  \| \partial_\omega \eta \|_4  \|\partial_\omega \eta\|_{17}+\| \partial_\omega\eta \|_6 \| \partial_\omega  \eta \|_{16}   )))
  &\mbox{for }i=16;\\
K^2m  \| \partial_\omega  \eta \|_{17}^2  + \|\eta\|_{18} (    \|\eta\|_{18}  (
\| u\|_3^2+ \| u\|_3^4  & \\
+m^4\| \partial_\omega  \eta \|_4^2 \|  \partial_\omega   \eta \|_6^2   ) +m^2 (\| \partial_\omega  \eta \|_3 \| \partial_\omega  \eta \|_6\|\eta\|_{17}\|u\|_3 &\\
 +  \| \partial_\omega \eta \|_6 \| \partial_\omega \eta \|_{17} \|u\|_3 ) )  + m^2 \| \partial_\omega \eta \|_4 \| \partial_\omega \eta \|_{17}^{1/2}\|\eta\|_{18}^{3/2}\|u\|_3 &\mbox{for }i=17.
\end{cases}
\label{qwessebdaiseqinM0846dfgssgsd}
 \end{align}
  \end{lem}
\begin{pf}
Let $i=0$, $\ldots$, $12$, $16$, $17$, and $\alpha$ satisfy $|\alpha|=i$.
 If $\alpha\neq 0$, then there exists a component, denoted by $\alpha_j$, such that $\alpha_j\neq0$, and thus we further define $\alpha^-$ and $\alpha^+$ as follows:
\begin{align}
\label{202004182208}
\alpha^-_k =\alpha^+_k=\alpha_k\mbox{ for }k\neq j,\ \alpha^-_j= \alpha_j-1
\mbox{ and } \alpha^+_j=\alpha_j+1.
\end{align}

Applying $ \partial^{ \alpha} $ to \eqref{s0106pnnnn}$_1$ yields
\begin{align}
\label{20221111620049}
 \partial^{ \alpha}(u_t-\nu\Delta  u-m^2\partial_\omega^2\eta
)=\partial^{ \alpha}(\mathcal{N}^\nu-\nabla_{\tilde{\ml{A}}}q - \nabla q).
\end{align}
Multiply the above identity by $\partial^{\alpha} u$ in $L^2$, and then using the integral by parts and \eqref{s0106pnnnn}$_2$,
we get
\begin{align}
& \frac{1}{2}\frac{\mm{d}}{\mm{d}t}  \|\partial^\alpha (u,m\partial_\omega  \partial^\alpha\eta)\|_0^2+ \nu\|\nabla \partial^\alpha u \|_0^2 = I_1+I_2, \label{qwe201808100154628xxx}
\end{align}
where $
I_1:= -  \int   \partial^{\alpha } \mathcal{N}^\nu_{j,l}  \partial_l \partial^{\alpha} u_j\mm{d}y$ and
  \begin{align}
  I_2:=-\int   \partial^\alpha \mm{div}_{\tilde{\mathcal{A}}} u \partial^\alpha  q   \mm{d}y
\begin{cases}
\displaystyle -\int    {u}\cdot  \nabla_{\tilde{\mathcal{A}}} q    \mm{d}y &\mbox{for }\alpha= 0;\\ \displaystyle
+\int   \partial^{\alpha^+} {u}\cdot\partial^{\alpha^-} \nabla_{\tilde{\mathcal{A}}} q    \mm{d}y&\mbox{for }\alpha\neq 0.
\end{cases}\nonumber
\end{align}

Exploiting H\"older's inequality, \eqref{202005sda08}, \eqref{200s12} and \eqref{apxxx1x}, we find that
\begin{align}
I_1\lesssim &\sum_{1\leqslant j, l\leqslant 3}
(\|  \partial_l u_j\|_{i}\| \mathcal{N}^\nu_{j,l}\|_i)\lesssim
\|u\|_{i+1} (1+\|\tilde{\mathcal{A}} \|_{2} )( \|\tilde{\mathcal{A}} \|_{2} \| u\|_{i+1}+\|\tilde{\mathcal{A}} \|_{i} \|u\|_3  ) \nonumber \\
\lesssim &
\| \eta\|_{i+1}\|u\|_3\|u\|_{i+1}  +\delta\|u\|_{i+1}^2  .\label{202sdfaf04101345}
\end{align}
Similarly, the integral $I_2$ can be estimated as follows.
\begin{align}
I_2\lesssim&\|\tilde{\mathcal{A}} \|_{i} \|  u\|_3 \|q\|_{i}+
 \|  u\|_{i+1} (\|\tilde{\mathcal{A}} \|_{\sigma(i-1)} \|  q\|_3  +\|\tilde{\mathcal{A}} \|_{2} \| q\|_{i}  )
 \nonumber \\
\lesssim&
 \| \eta\|_{i+1 } \|u\|_3\|q\|_{i}+ \|  u\|_{i+1}(\|  \eta\|_{\sigma(i-1)+1} \| q\|_3  +\| \eta \|_3 \| q\|_{i}  ) ,\label{20200614145dfs2}
   \end{align}
   see \eqref{2022111161031} for the definition of $\sigma$.

(1) \emph{Case $0\leqslant i\leqslant 3$}.

Exploiting \eqref{2020211013200910} and the interpolation inequality \eqref{202108261405}, we further derive from \eqref{202sdfaf04101345} with $0\leqslant i\leqslant 3$ that
\begin{align}
I_1\lesssim
\|\partial_\omega \eta\|_{i}^{1/2}\|\partial_\omega \eta\|_{8+i}^{1/2}\|u\|_{i+1}^{3/2} \|u\|_{5-i}^{1/2} +\delta\|u\|_{i+1}^2  .\label{20201345}
\end{align}

Similarly, thanks to \eqref{2020211013200910},  \eqref{2017020614181721} with $i=1$, we can further derive from  \eqref{20200614145dfs2} with $0\leqslant i\leqslant 3$ that
\begin{align}
I_2\lesssim&
 (\|\eta\|_{3} \| u\|_{i+1}+\|  \eta\|_{i+1} \| u\|_3 )( \|  u\|_2 \|  u\|_3 \nonumber \\
 &+ m^2  (\| \partial_\omega  \eta \|_2 \|  \partial_\omega \eta \|_3
    + \|  \eta \|_3 \| \partial_\omega^2 \eta \|_2))\nonumber \\
    \lesssim&   \|  \eta\|_{4} \| u\|_4^3+
   m^2 \|\partial_\omega \eta\|_3\|\partial_\omega \eta\|_6 \| \partial_\omega \eta\|_7\| u\|_{4}
\nonumber \\
\lesssim& \|\eta\|_{4}\|u\|_{10-2i}\| u\|_{i+1}^2+  m^2 \| \partial_\omega \eta \|_i^2  \| \partial_\omega \eta\|_{16-2i}\| u\|_{4} . \label{202021010131631}
   \end{align}
Putting  the above two estimates into \eqref{qwe201808100154628xxx} for $|\alpha|=i$, and then using  Young inequality, \eqref{20160614fdsa19asfda57x}  and \eqref{aprpxx1x}, we get \eqref{qwessebdaiseqinM0846dfgssgsd} with $0\leqslant i\leqslant 3$ for sufficiently small $\delta$.

 (2) \emph{Case $4\leqslant i\leqslant 12$}.

We can derive from  \eqref{20200614145dfs2} with $4 \leqslant i\leqslant 12$ that
\begin{align}
I_2\lesssim&
 \| \eta\|_{i}\|  u\|_{i+1} ( \|    u\|_2 \| u\|_3  + m^2  (\| \partial_\omega  \eta \|_2 \|  \partial_\omega \eta \|_3
    + \|  \eta \|_3 \| \partial_\omega^2\eta \|_2)) \nonumber \\
  &+(\|\eta \|_3 \|u\|_{i+1}+
 \| \eta\|_{i+1}\|u\|_3  )(\|   u\|_2 (\| u\|_i+ \|\eta\|_i\| u\|_3)+ \|  u\|_3 \|   u\|_{i-1}
\nonumber\\
& + m^2 (\|  \eta \|_3\| \partial_\omega^2 \eta \|_{i-1} + \| \partial_\omega  \eta \|_3 \|  \partial_\omega \eta \|_{i-1} + \|  \eta \|_{i-1}(\| \partial_\omega^2 \eta \|_3+  \| \partial_\omega  \eta \|_3^2 ) \nonumber \\
&+  \| \eta\|_{i}  ( \| \partial_\omega   \eta \|_2 \|  \partial_\omega \eta \|_3+ \| \eta \|_3 \| \partial_\omega^2 \eta \|_2
    ))) \nonumber \\
\lesssim&\| \eta\|_{i+1} \|  u\|_{i+1}(1+\|\eta\|_{i})(\|    u\|_3 \| u\|_{i}  + m^2 \| \partial_\omega\eta \|_6   \|  \partial_\omega \eta \|_{i} )
   \nonumber \\
   &   +m^2\| \partial_\omega  \eta \|_4 \| \partial_\omega \eta \|_{i+2}\| \partial_\omega \eta\|_{i+4}\|u\|_3
 ,
   \end{align}
   where we have used \eqref{2017020614181721} with $i=1$, $10$, in the first inequality,  and \eqref{2020211013200910} in the   last inequality.

 Putting  \eqref{202sdfaf04101345} with  $4\leqslant i\leqslant 12$ and  the above estimate into \eqref{qwe201808100154628xxx} for $|\alpha|=i$, and then using  \eqref{20160614fdsa19asfda57x},   \eqref{aprpiosesnewxxxx1x}, \eqref{aprpxx1x}  and Young's inequality, we get \eqref{qwessebdaiseqinM0846dfgssgsd} with $4\leqslant i\leqslant 12$.

(3) \emph{Case $i=16$}.

We can derive from  \eqref{20200614145dfs2} with $i=16$ that
 \begin{align}
I_2\lesssim&
 \|\eta\|_{16}\|u\|_{17}( \|    u\|_2 \|  u\|_3  + m^2  (\| \partial_\omega  \eta \|_2 \|  \partial_\omega   \eta \|_3
    + \|  \eta \|_3  \| \partial_\omega^2 \eta \|_2)) \nonumber \\
  &+(\| \eta \|_3 \|u\|_{17}+
 \|\eta\|_{17}\|u\|_3 )(\|   u\|_2 (\|u\|_{16}+ \|\eta\|_{16}\| u\|_3)+ \|  u\|_3 \|   u\|_{15}
\nonumber\\
& + m^2 (  \| \eta \|_3 \| \partial_\omega^2 \eta \|_{15} + \| \partial_\omega  \eta \|_3 \|  \partial_\omega  \eta \|_{15}   + \| \eta \|_{15}(\| \partial_\omega^2 \eta \|_3+\| \partial_\omega  \eta \|_3^2 ) \nonumber  \\
&+   \| \eta\|_{16}   (\| \partial_\omega  \eta \|_2 \|  \partial_\omega   \eta \|_3
    + \|  \eta \|_3  \| \partial_\omega^2 \eta \|_2)))\nonumber\\
\lesssim& \|u\|_{17}( \|\eta\|_{17} (1+ \| \eta \|_3 )( \| u\|_3^2+ m^2    \| \partial_\omega \eta \|_4   \|  \partial_\omega   \eta \|_6 )\nonumber \\
& +
\| \eta \|_3 ( \|   u\|_3  \|u\|_{16}+  m^2\| \partial_\omega\eta \|_6 \| \partial_\omega \eta \|_{16} ))   +\|\eta\|_{17}(\|\eta\|_{16}\|   u\|_3( \|   u\|_2 \|   u\|_3 \nonumber \\
  &  +m^2\| \partial_\omega  \eta \|_3 \|\partial_\omega \eta \|_6)+ m^2\| \partial_\omega\eta \|_6 \| \partial_\omega  \eta \|_{16} \|u\|_3  )+    m^2\| \partial_\omega \eta \|_4\|\eta\|_{15} \|\eta\|_{17}\|   u\|_3\nonumber\\
\lesssim&
\delta  (\|u\|_{17}^2+ m \| \partial_\omega \eta \|_{16}\|u\|_{17} )+ \|\eta\|_{17} ( \|u\|_{17} ( \| u\|_3^2+ m^2    \| \partial_\omega \eta \|_4   \|  \partial_\omega   \eta \|_6 ) \nonumber \\
& + \|\eta\|_{16} \|   u\|_3(  \|   u\|_2 \|   u\|_3   +m^2\| \partial_\omega  \eta \|_3 \|\partial_\omega \eta \|_6) + m^2\| \partial_\omega\eta \|_6 \| \partial_\omega  \eta \|_{16} \|u\|_3  )\nonumber \\
  &+   m^2 \| \partial_\omega \eta \|_4  \|\partial_\omega \eta\|_{17}\|\eta\|_{18} \|   u\|_3
 ,
   \end{align}
   where we have used \eqref{2017020614181721} with $i=1$, $14$, in the first inequality,  \eqref{2020211013200910} in the second inequality, and
 \eqref{2020211013200910}, \eqref{aprpiosesnewxxxx1x}, \eqref{aprpxx1x},  the interpolation inequality \eqref{202108261405}  in the last inequality.

 Putting  \eqref{202sdfaf04101345} with  $i=16$ and  the above estimate into \eqref{qwe201808100154628xxx} for $|\alpha|=16$, and then using  \eqref{202021100321520}, \eqref{20160614fdsa19asfda57x}  and Young's inequality, we get \eqref{qwessebdaiseqinM0846dfgssgsd} with $i=16$.

(4) \emph{Case $i=17$}.

We can derive from  \eqref{20200614145dfs2} with $i=17$ that
 \begin{align}
I_2\lesssim&
 \|\eta\|_{17}\|u\|_{18}( \|    u\|_2 \|  u\|_3  + m^2  (\| \partial_\omega  \eta \|_2 \|  \partial_\omega   \eta \|_3
    + \|  \eta \|_3  \| \partial_\omega^2 \eta \|_2)) \nonumber \\
  &+(\| \eta \|_3 \|u\|_{18}+
 \|\eta\|_{18}\|u\|_3 )((\|   u\|_2 (\|u\|_{17}+ \|\eta\|_{17}\| u\|_3)+ \|  u\|_3 \|   u\|_{16}
\nonumber\\
& + m^2 ( \| \eta \|_3 \| \partial_\omega^2 \eta \|_{16} + \| \partial_\omega  \eta \|_3 \|  \partial_\omega  \eta \|_{16}    + \| \eta\|_{16}( \| \partial_\omega^2 \eta \|_3+  \| \partial_\omega  \eta \|_3^2  )
 \nonumber  \\
 &+  \| \eta\|_{17}   (\| \partial_\omega  \eta \|_2 \|  \partial_\omega   \eta \|_3
    + \|  \eta \|_3  \| \partial_\omega^2 \eta \|_2)))\nonumber \\
\lesssim&\|u\|_{18}( \|\eta\|_{18} (1+ \| \eta \|_3  ) (  \|  u\|_3^2
 + m^2\| \partial_\omega  \eta \|_4 \|  \partial_\omega   \eta \|_6)  \nonumber \\
 &+
\| \eta \|_3 (\|   u\|_3  \|u\|_{17} +  m^2\| \partial_\omega\eta \|_6 \| \partial_\omega \eta \|_{17}))   +\|\eta\|_{18}(\|\eta\|_{17}\|   u\|_3( \|   u\|_2 \|   u\|_3 \nonumber \\
  &  +m^2\| \partial_\omega  \eta \|_3 \|\partial_\omega \eta \|_6)+ m^2\| \partial_\omega\eta \|_6 \| \partial_\omega  \eta \|_{17} \|u\|_3  )   +  m^2 \| \partial_\omega \eta \|_4 \| \eta \|_{16}\|\eta\|_{18}\|u\|_3 \nonumber
  \\ \lesssim& \delta\|u\|_{18}^2 +K\sqrt{m}   \| \partial_\omega  \eta \|_{17}\|u\|_{18}   + \|\eta\|_{18} (\|u\|_{18}   (  \|  u\|_3^2
 + m^2\| \partial_\omega  \eta \|_4 \|  \partial_\omega   \eta \|_6)  \nonumber \\
 & +\|\eta\|_{17}\|u\|_3 ( \|   u\|_2 \|u\|_3+m^2   \| \partial_\omega  \eta \|_3 \| \partial_\omega  \eta \|_6   ) +m^2 \| \partial_\omega \eta \|_6 \| \partial_\omega \eta \|_{17} \|u\|_3 ) \nonumber \\
  & + m^2 \| \partial_\omega \eta \|_4 \| \partial_\omega \eta \|_{17}^{1/2}\|\eta\|_{18}^{3/2}\|u\|_3,  \label{2022111040851}
   \end{align}
    where we have used \eqref{2017020614181721} with $i=1$, $15$, in the first inequality,  \eqref{2020211013200910} in the second inequality and
  \eqref{aprpiosesnewxxxx1x}, \eqref{aprpxx1x} the interpolation inequality in the last inequality.

Putting \eqref{202sdfaf04101345} with  $i=17$  and  the above estimate into \eqref{qwe201808100154628xxx} for $|\alpha|=17$, and then using \eqref{20160614fdsa19asfda57x}  and Young's inequality, we get \eqref{qwessebdaiseqinM0846dfgssgsd} with $i=17$.
 \hfill$\Box$
\end{pf}

\begin{lem}
\label{qwe201612132242nxsfssdfsxx}
 Under the conditions \eqref{aprpiosesnewxxxx}--\eqref{aprpiose1snewxxxxz} with sufficiently small $\delta$, we have \begin{align}
 & \frac{\mm{d}}{\mm{d}t}\left(\frac{\nu}{2}\|\nabla^{{i+1}}\eta\|_0^2  +   \sum_{|\alpha|=i} \int \partial^\alpha \eta\cdot \partial^\alpha u\mm{d}y \right) +  \| m \partial_\omega \eta \|_i^2 \nonumber \\
& \leqslant
   2\|   u \|_{ i}^2+\times \begin{cases}
  c\delta  \|   u\|_{i+1}^2 \qquad \qquad  \qquad \quad\quad   \mbox{\ for }0\leqslant i\leqslant 12;&\\
c\| \eta\|_{i+1} (\| \eta\|_{i+1}(m^2 \|
 \partial_\omega \eta \|_6^2&\\
  +\| u\|_3(1 + \| u\|_3))&\\
 + \begin{cases}
 \| \eta \|_3 \| u\|_{i+1}  )  &\mbox{for }i=16;\\
 \| \eta \|_3^{2/3} \|\partial_\omega \eta \|_6^{1/3} \| u\|_{i+1}  )  &\mbox{for }i=17.
  \end{cases}&
\end{cases}
  \label{qweLem:03xx}
 \end{align}
\end{lem}
\begin{pf}
Let $i=|{\alpha}|= 0$, $\ldots$, $12$, $16$ and $17$.
Multiplying \eqref{20221111620049} by $\partial^{\alpha}\eta$ in $L^2$,
we get that
\begin{align}
 &\frac{\mm{d}}{\mm{d}t}  \left(\frac{\nu}{2}\|\nabla \partial^\alpha \eta \|_0^2+ \int \partial^\alpha \eta \cdot\partial^\alpha u\mm{d}y \right)
  +  \|m\partial_\omega\partial^\alpha \eta \|_0^2 = \|\partial^\alpha u\|_0^2+I_3+I_4, \label{qwe2010154628}
\end{align}
where
$I_3:=-\int  \partial^\alpha {\mathcal{N}}^\nu_{j,l}\cdot\partial_l\partial^\alpha\eta_j\mm{d}y$
 and
  \begin{align}
I_4:=  \int \partial^\alpha q\partial^\alpha (r^\eta_2+r^\eta_3) \mm{d}y
\begin{cases}
\displaystyle -\int  \eta\cdot  \nabla_{\tilde{\mathcal{A}}} q    \mm{d}y &\mbox{for }\alpha= 0;\\ \displaystyle
+\int   \partial^{\alpha^+} {\eta}\cdot\partial^{\alpha^-} \nabla_{\tilde{\mathcal{A}}} q    \mm{d}y&\mbox{for }\alpha\neq 0.
\end{cases}\nonumber
\end{align}
with $\alpha^-$ and $\alpha^+$ being defined by \eqref{202004182208}.

Following the arguments of \eqref{202sdfaf04101345} and \eqref{20200614145dfs2} with slight modifications,
  we find that
\begin{align}
 I_3 \lesssim  &  \|\eta \|_3\|\eta\|_{i+1} \|u\|_{i+1}
 +\begin{cases}
0&\mbox{for }i=0,\ 1;\\
 \|\eta\|_{i+1}^2\|u\|_3&\mbox{for }i\geqslant 2 ,
  \end{cases} \label{202005safas0119091}\\
 I_4 \lesssim &\|\eta\|_{i+1}  (\|\eta \|_3 \|q\|_{i} + \| \eta\|_{\sigma(i-1)+1 }  \|q\|_{3}). \label{2020050119091}
\end{align}

(1) \emph{Case $0\leqslant i\leqslant 3$}.

 Exploiting  \eqref{2020211013200910}, \eqref{202108261405} and \eqref{2017020614181721}, we derive from \eqref{202005safas0119091} and   \eqref{2020050119091} with $i=3$ that
\begin{align}
 I_3  \lesssim  &\|u\|_{i+1} (\|\partial_\omega\eta\|_{6}\|\partial_\omega\eta\|_{i+4} +\| \partial_\omega\eta\|_{i+4}^2)   \nonumber
 \\
 \lesssim & \|\partial_\omega\eta\|_{i}\|u\|_{i+1} (\|\partial_\omega\eta\|_{10}+ \| \partial_\omega\eta\|_{i+8})\label{202110211607}
\end{align}
and
\begin{align}
I_4\lesssim&
 \| \eta\|_3 \| \eta\|_4( \|u\|_2 \| u\|_3 +m^2  ( \| \partial_\omega   \eta \|_2 \|  \partial_\omega \eta \|_3+\|  \eta \|_3 \| \partial_\omega^2 \eta \|_2 ))  \nonumber \\
 \lesssim&  \| \partial_\omega\eta\|_7^2(  \| u\|_3^2 +
 m^2  \|  \partial_\omega \eta \|_6 \| \partial_\omega \eta \|_3 )    \nonumber\\
 \lesssim& \| \partial_\omega  \eta\|_{i}^{2(3-i)/(10-2i)} \| \partial_\omega  \eta\|_{10-i}^{2(7-i)/(10-2i)}\|  u\|_{i}^{2(7-i)/(10-2i)} \|  u\|_{10-i}^{2(3-i)/(10-2i)}   \nonumber\\& + m^2\|\partial_\omega\eta\|_{i}^{25/(12-i)}\| \partial_\omega\eta\|_{12}^{(23-4i)/(12-i) }, \nonumber
   \end{align}
   where $0< 2(3-i)(10-2i)^{-1}\leqslant 1$ for $0\leqslant i\leqslant 2$ and $1< (23-4i)/(12-i)< 2$ for any $0\leqslant i\leqslant 3$.
Inserting the above two estimates into \eqref{qwe2010154628} with $i=3$, and then using  \eqref{apxxx1x}  and Young's inequality, we obtain \eqref{qweLem:03xx} with $0\leqslant i\leqslant 3$.

(2) \emph{Case $4\leqslant i\leqslant 12$.}

Similarly to \eqref{202110211607}, we have, for $4\leqslant i\leqslant 12$,
\begin{align}
 I_3 \lesssim  & \| \eta\|_{i+1}^2   \|u\|_{i+1}  \lesssim
  \| \partial_\omega\eta\|_{i}\|  \eta\|_{i+5}    \|u\|_{i+1} \lesssim \delta m \| \partial_\omega\eta\|_{i}   \|u\|_{i+1}, \label{20200091}
\end{align}
where we have used \eqref{apxxx1x} in the last inequality.

Noting that
\begin{align}
\|\eta\|_i\|\partial_\omega \eta\|_6\lesssim \|\partial_\omega \eta\|_i\|\partial_\omega \eta\|_{9}\mbox{ for }i=4,\ 5, \label{2022111101948}
\end{align}
thus making use of \eqref{2020211013200910}, \eqref{2017020614181721} with $i=1$, $10$,  \eqref{aprpiosesnewxxxx1x}, \eqref{aprpxx1x} and \eqref{2022111101948},  we can derive from \eqref{2020050119091} that,
for   $4\leqslant i\leqslant 12$,
\begin{align}
|I_4| \lesssim  & \|\eta \|_3 \|\eta\|_{i+1} (\|   u\|_2 (\| u\|_{i}+ \|\eta\|_{i}\| u\|_3)+ \|  u\|_3 \|   u\|_{i-1}
\nonumber\\
& + m^2 ( \|  \eta \|_3\| \partial_\omega^2 \eta \|_{i-1} + \| \partial_\omega  \eta \|_3 \|  \partial_\omega \eta \|_{i-1}     +\| \eta\|_{i-1}  (\| \partial_\omega^2 \eta \|_3+  \| \partial_\omega  \eta \|_3^2 )   \nonumber \\
&+   \| \eta\|_{i} (\| \partial_\omega   \eta \|_2 \|  \partial_\omega \eta \|_3
    + \| \eta \|_3 \| \partial_\omega^2 \eta \|_2))  )\nonumber \\
 &+    \| \eta\|_{i}\|  \eta\|_{i+1} ( \|    u\|_2 \| u\|_3  + m^2 (\| \partial_\omega  \eta \|_2 \|  \partial_\omega   \eta \|_3
    + \|  \eta \|_3  \| \partial_\omega^2 \eta \|_2)) \nonumber \\
\lesssim&   \|\eta \|_3 \|\eta\|_{i+1} ( \|  u\|_3 \|   u\|_{i}+m^2 \| \partial_\omega \eta \|_6\| \partial_\omega \eta \|_{i})\nonumber \\
  &+ \| \eta\|_{i}\|  \eta\|_{i+1}(1+\|\eta \|_3 ) ( \|    u\|_2 \| u\|_3  + m^2\| \partial_\omega \eta \|_4  \| \partial_\omega \eta \|_6 ) \nonumber \\
\lesssim& \delta ( \| u\|_{i}^2+ m^2  \| \partial_\omega \eta \|_{i}^2   )  \label{202111112111}
   \end{align}
Inserting \eqref{20200091} and \eqref{202111112111}  into \eqref{qwe2010154628} with $i=12$, and then using    Young's inequality, we obtain \eqref{qweLem:03xx} with $4\leqslant i\leqslant 12$.

(3) \emph{Case $i=16$}.

We can derive from \eqref{2020050119091} with $i=16$ that
\begin{align}
|I_4| \lesssim  &\| \eta \|_3 \|\eta\|_{17} (\|   u\|_2 (\|u\|_{16}+ \|\eta\|_{16}\| u\|_3)+ \|  u\|_3 \|   u\|_{15}
\nonumber\\
& + m^2 (\| \eta \|_3 \| \partial_\omega^2 \eta \|_{15} + \| \partial_\omega  \eta \|_3 \|  \partial_\omega  \eta \|_{15}    +\| \eta\|_{15} (
   \| \partial_\omega^2 \eta \|_3 +\| \partial_\omega  \eta \|_3^2)  \nonumber  \\
&+  \| \eta\|_{16}   (\| \partial_\omega  \eta \|_2 \|  \partial_\omega   \eta \|_3
    + \|  \eta \|_3  \| \partial_\omega^2 \eta \|_2)))\nonumber \\
  &+
  \|\eta\|_{16}\|\eta\|_{17}( \|    u\|_2 \|  u\|_3  + m^2(\| \partial_\omega  \eta \|_2 \|  \partial_\omega   \eta \|_3
    + \|  \eta \|_3  \| \partial_\omega^2 \eta \|_2)) \nonumber \\
\lesssim& \| \eta \|_3 \|\eta\|_{17} (\|   u\|_3  \|u\|_{16}
 + m^2  \| \partial_\omega \eta \|_6 \| \partial_\omega  \eta \|_{16}   ) \nonumber \\
  & +\|\eta\|_{16}\|\eta\|_{17}(1 +\| \eta \|_3)( \|    u\|_2 \|  u\|_3  + m^2   \| \partial_\omega \eta \|_4\|\partial_\omega  \eta \|_6
   )\nonumber  \\
\lesssim&  \|\eta\|_{17}^2(   \|  u\|_3 ^2 + m^2   \|\partial_\omega  \eta \|_6^2
   ) +  \delta \|(u,m \partial_\omega  \eta) \|_{16} ^2 , \label{202110181927}
   \end{align}
   where we have used    \eqref{2017020614181721} with $i=1$, $14$ in the first inequality,   \eqref{2020211013200910} in the second inequality
  and
  \eqref{aprpiosesnewxxxx1x} and \eqref{aprpxx1x} in the last inequality.
Inserting \eqref{202005safas0119091} with $i=16$ and the above estimate into \eqref{qwe2010154628} for $|\alpha|=16$ yields  \eqref{qweLem:03xx} with $i=16$.

(4) \emph{Case $i=17$}.

We can derive from \eqref{2020050119091} with $i=17$ that
\begin{align}
|I_4| \lesssim
  &  \| \eta \|_3 \|\eta\|_{18} (\|   u\|_2 (\|u\|_{17}+ \|\eta\|_{17}\| u\|_3)+ \|  u\|_3 \|   u\|_{16}
\nonumber\\
& + m^2 (  \| \eta \|_3 \| \partial_\omega^2 \eta \|_{16} + \| \partial_\omega  \eta \|_3 \|  \partial_\omega  \eta \|_{16}    +\| \eta\|_{16}(
\| \partial_\omega^2\eta \|_3+  \| \partial_\omega  \eta \|_3^2) \nonumber  \\
&   +   \| \eta\|_{17}  (\| \partial_\omega  \eta \|_2 \|  \partial_\omega   \eta \|_3
    + \|  \eta \|_3  \| \partial_\omega^2 \eta \|_2) ))\nonumber \\
&+
 \|\eta\|_{17}\|\eta\|_{18}( \|    u\|_2 \|  u\|_3  + m^2 (\| \partial_\omega  \eta \|_2 \|  \partial_\omega   \eta \|_3
    + \|  \eta \|_3  \| \partial_\omega^2 \eta \|_2)) \nonumber \\
\lesssim&  \| \eta \|_3 \|\eta\|_{18} (\|   u\|_3  \|u\|_{17}  + m^2   \| \partial_\omega\eta \|_6\| \partial_\omega \eta \|_{17}    ) \nonumber \\
  &+\|\eta\|_{17}\|\eta\|_{18}(1+\| \eta \|_3)( \|    u\|_2 \|  u\|_3  + m^2    \| \partial_\omega \eta \|_4\|
 \partial_\omega \eta \|_6) \nonumber \\
\lesssim&  \|\eta\|_{18}^2 ( \|  u\|_3^2  + m^2   \|
 \partial_\omega \eta \|_6^2)  +   \delta \|(u, m \partial_\omega \eta )\|_{17}^2      , \label{2021127}
   \end{align}
   where we have used  \eqref{2017020614181721} with $i=1$, $15$ in the first inequality, \eqref{2020211013200910}  in the second inequality,   and
  \eqref{aprpiosesnewxxxx1x} and \eqref{aprpxx1x} in the last inequality.
Inserting \eqref{202005safas0119091} with $i=17$ and the above estimate into \eqref{qwe2010154628} for $|\alpha|=17$, and then using \eqref{2020211013200910}, we  get  \eqref{qweLem:03xx} with $i=17$. \hfill$\Box$
 \end{pf}

\subsection{Stability estimates} \label{202006141016}
With the energy estimates in Lemmas \ref{qwe201612132242nn}--\ref{qwe201612132242nxsfssdfsxx} in hand, we are in the position to establish
the \emph{a priori} estimate \eqref{aprpiose1}.

To begin with, we  use \eqref{qwessebdaiseqinM0846dfgssgsd} and \eqref{qweLem:03xx} with $0\leqslant i\leqslant 12$ to build the following $i$-th layer  energy inequality:
\begin{align}
&\frac{\mm{d}}{\mm{d}t}    \tilde{\mathcal{E}}_{i}+
c \mathcal{D}_{i}   \lesssim  0
\label{20191sadf1262012}
\end{align}
where
\begin{align}
\tilde{\mathcal{E}}_{i}  :=c \|\nabla^{i}(u,m\partial_\omega  \eta)\|_0^2 +\frac{\nu}{2}\|\nabla^{i+1} \eta\|_0^2 +
\sum_{  |\alpha| = i} \int  \partial^\alpha \eta \cdot \partial^\alpha u\mm{d}y
\label{202211052115}
  \end{align}
and ${\mathcal{E}}_{i}$ satisfies
\begin{align}
 \mathcal{E}_{i}\lesssim   \tilde{\mathcal{E}}_{i}\lesssim  \mathcal{E}_{i}.  \label{202006141102}
\end{align}
Integrating \eqref{20191sadf1262012}   over $(0,t)$ yields
\begin{align}
& \mathcal{E}_{i} +
  \int_0^t\mathcal{D}_{i} \mm{d}\tau\lesssim \mathcal{E}^{0}_{i}.
\label{2019112safdsafsaf6201sdaf2}
\end{align}

By \eqref{2020211013200910}, we see that, for $0\leqslant i\leqslant 12$,
\begin{align}
\mathcal{E}_{i} \lesssim  (1+m^{-2})\mathcal{D}_{i+4}.
\label{202002111102021}
\end{align}
 Thus  we further derive the following lower-order energy inequality from \eqref{20191sadf1262012}, \eqref{202006141102} and \eqref{202002111102021}:
\begin{align}
& \sum_{i=0}^{3} \left(  d_{i }(1+m^{-2})^i\frac{\mm{d}}{\mm{d}t}\langle t\rangle^{ (3-i)}\tilde{\mathcal{E}}_{4i}+
h_{i}(1+m^{-2})^i\langle t\rangle^{ (3-i)} \mathcal{D}_{4i}   \right)\lesssim  0
\label{20112}
\end{align}
for some constants $d_i$  and $h_i$ depending on $\nu$ and $\omega$.
Integrating \eqref{20112}  over $(0,t)$, and then using \eqref{202006141102}, we get
\begin{align}
& \sum_{i=0}^{3} \left((1+m^{-2})^i \langle t\rangle^{ (3-i)} {\mathcal{E}}_{4i}+ (1+m^{-2})^i\int_0^t\langle t\rangle^{ (3-i)}\mathcal{D}_{4i}  \mm{d}\tau\right)\lesssim \Xi,
\label{20191201sdaf2}
\end{align}
where $\Xi$ is defined by \eqref{20202111110854}. Using the interpolation inequality, we can derive from \eqref{2019112safdsafsaf6201sdaf2} and \eqref{20191201sdaf2} that, for $0\leqslant i\leqslant 12$,
\begin{align}
\label{2x62012}
\mathcal{E}_{i}\lesssim\left(\mathcal{E}^{0}_{12}\right)^{i/12} \left( \langle t\rangle^{ -3} \Xi\right)^{(12-i)/12}  .
\end{align}

It is easy see from \eqref{qwessebdaiseqinM0846dfgssgsd} and \eqref{qweLem:03xx} that
\begin{align}
&\frac{\mm{d}}{\mm{d}t} \|\nabla^{i} (u, m\partial_\omega \eta ) \|_0^2   + \nu\| u \|_{i+1}^2
\nonumber \\
&\lesssim \begin{cases}
 \delta \mathcal{D}_H + \mathcal{E}_H\mathcal{D}_6(1+\sqrt{\mathcal{E}_3}+\mathcal{E}_4)
 +\sqrt{\mathcal{E}_3\mathcal{E}_4\mathcal{E}_H\mathcal{D}_H}
  &\mbox{for }i=16;\\
K^2m^{-1} \mathcal{D}_H
+\mathcal{E}_H (\mathcal{D}_6( m^{2/3} +K+ K \sqrt{\mathcal{E}}_4 )+ K\sqrt{\mathcal{E}}_3) &\mbox{for }i=17.
\end{cases}
\label{qwesgsd}
 \end{align}
 and  \begin{align}
 & \frac{\mm{d}}{\mm{d}t}\left(\frac{\nu}{2}\|\nabla^{{i+1}}\eta\|_0^2  +   \sum_{|\alpha|=i} \int \partial^\alpha \eta\cdot \partial^\alpha u\mm{d}y \right) +  \| m \partial_\omega \eta \|_i^2 \nonumber \\
& \leqslant
   2\|   u \|_{ i}^2+c ( \sqrt{\mathcal{E}_H\mathcal{D}_H} (\sqrt{\mathcal{E}_3} +\mathcal{E}_3^{1/3}\mathcal{D}_6^{1/6} )+ \mathcal{E}_H(\sqrt{\mathcal{E}_3}+\mathcal{D}_6))
  \label{qweLemxx}\mbox{ for }i=16,\ 17,
 \end{align}

Noting that, by \eqref{202021100321520} and \eqref{aprpiose1snewxxxxz},
\begin{align}
\frac{K^{2}}{m^{5/3}} \lesssim_0 \delta^{5/3}\mbox{ and } \frac{K}{m^{2/3}} \lesssim_0 \delta^{2/3}. \label{20221110510106}
\end{align}
Exploiting \eqref{202021100321520}, Poinc\'are's inequality \eqref{20160614fdsa19asfda57x}, \eqref{20221110510106} and Young's inequality, we derive the following highest-order energy inequality from \eqref{qwesgsd} and \eqref{qweLemxx}:
\begin{align}
 \frac{\mm{d}}{\mm{d}t}    \tilde{\mathcal{E}}_H+
c \mathcal{D}_{H}  \lesssim     \mathcal{E}_H(  \sqrt{\mathcal{E}_{3} }
 +\mathcal{E}_{3}\mathcal{E}_4+ (1+\mathcal{E}_{4} ) \mathcal{D}_{6}),
\label{20191sa012}
\end{align}
where
$$
\begin{aligned}
\tilde{\mathcal{E}}_H:= &c
 \|\nabla^{16}  (u,m\partial_\omega  \eta)\|_0^2+m^{-2/3}\|\nabla^{17}  (u,m\partial_\omega  \eta)\|_0^2\nonumber \\
 &+ \frac{\nu}{2}\|\nabla^{17}( \eta, \nabla \eta )\|_0^2+
\sum_{16\leqslant   |\alpha|\leqslant 17} \int  \partial^\alpha  \eta  \cdot \partial^\alpha  u\mm{d}y,
\end{aligned}$$
satisfying
\begin{align}
\mathcal{E}_H\lesssim   \tilde{\mathcal{E}}_H \lesssim \mathcal{E}_H.  \label{202102}
\end{align}

Applying Gronwall's lemma to \eqref{20191sa012}, and then using   \eqref{2019112safdsafsaf6201sdaf2}, \eqref{2x62012} with $i=3$ and \eqref{202102}, we arrive at that there exists a constant $\delta_1\in(0,1]$ such that for any $\delta\leqslant \delta_1$,
\begin{align}
 \mathcal{E}_{H}  \lesssim   \mathcal{E}_{H}^0 e^{c \int_0^T (  \sqrt{\mathcal{E}_{3} }
 +\mathcal{E}_{3}\mathcal{E}_4+ (1+\mathcal{E}_{4} ) \mathcal{D}_{6})\mm{d}\tau}
 \leqslant  &c_1\mathcal{E}_{H}^0 e^{c_2\vartheta}/4\mbox{ for any }0\leqslant t\leqslant T,
\label{20202110261936}
\end{align}
where $c_1\geqslant 4$ and $\vartheta $ is defined by \eqref{20211051042}. In addition, thanks to \eqref{20202110261936}, we further derive from \eqref{20191sa012} that
\begin{align}
 \int_0^t \mathcal{D}_{H} \mm{d}\tau
\lesssim  \mathcal{E}_{H}^0 (1+e^{ {c}_2\vartheta  } \vartheta  ).\label{20190sfda5sadfa041053}
\end{align} Now we take
\begin{equation}
\label{201911262060}
K:=  \sqrt{ c_1   \mathcal{E}_{H}^0 e^{ {c}_2\vartheta } }>0,
\end{equation}
we immediately obtain the desired \emph{a priori} stability estimate \eqref{aprpiose1} from \eqref{20202110261936}
 under the \emph{a priori} assumption  \eqref{aprpiosesnewxxxx} with the relative smallness condition \eqref{aprpiose1snewxxxxz}
for any  $\delta\leqslant \delta_1$.

%Finally, we further derive the quicker decay-in-time \eqref{20190safd5041053}. By the interpolation inequality and \eqref{2019112safdsafsaf6201sdaf2}, $$\|\eta\|_{4}\lesssim \|\partial_\omega \eta\|_{3}^{7/9}\|\eta\|_{18}^{2/9}\lesssim \|\partial_\omega \eta\|_{3}^{7/9}K^{2/9}. $$Thus $$  {\mathcal{E}}_1 \lesssim (m^{-18/13}(({\mathcal{E}}_{12}^0))^{4/13}+ (\mathcal{E}^0_3)^{4/13} )\mathcal{D}_3 )^{ 9/13},$$ which yields$$ ( {\mathcal{E}}_{L})^{ 13/9} / (m^{-18/13}(({\mathcal{E}}_{12}^0))^{4/13}+ (\mathcal{E}^0_3)^{4/13})\lesssim (m^{-18/13}+ (\mathcal{E}^0_3)^{4/13})\mathcal{D}_3. $$Putting the above estimate into  \eqref{20191sadf1262012} with $i=1$ yields that\begin{align}&\frac{\mm{d}}{\mm{d}t}    {\mathcal{E}}_1+ ( {\mathcal{E}}_1)^{ 13/9} / (m^{-18/13}(({\mathcal{E}}_{12}^0))^{4/13}+ (\mathcal{E}^0_3)^{4/13}) \lesssim  0. \nonumber\end{align}Applying the decay estimate \eqref{20201sdfa0261646} to the above inequality, we arrive at\begin{align}{\mathcal{E}}_1\lesssim_1 {\mathcal{E}}_1 \lesssim   (({\mathcal{E}}^0_1)^{-9/4}+ 9 t/4 (m^{-18/13}\sup_{t\in I_T} (({\mathcal{E}}_{12}^0))^{4/13}+ (\mathcal{E}^0_3)^{4/13}) )^{- 9 /4} .\label{2x62012}\end{align}
\subsection{Proof of Theorem \ref{201904301948}} \label{20206141015}
We start with introducing a local (-in-time) well-posedness result for the initial value problem \eqref{01dsaf16asdfasf}--\eqref{01dsaf16asdfasfsaf} and a result concerning diffeomorphism mappings.
\begin{pro} \label{pro:0401nxdxx}
Let $(\eta^0,u^0)\in H^{18}\times H^{17}$ satisfy $\|(\nabla\eta^0,u^0)\|_{17}\leqslant B$ and
$\mm{div}_{\mathcal{A}^0}u^0=0$, where $B$ is a positive constant, $\zeta^0:= \eta^0+y$ and $\mathcal{A}^0$ is defined by $\zeta^0$.
Then there is a constant $\delta_2\in (0,1]$, such that for any $(\eta^0,u^0)$ satisfying
\begin{align}
\|\nabla \eta^0\|_{12} \leqslant \delta_2, \label{201912261426}
\end{align}
 there exist a local existence time $T>0$ (depending possibly on $B$, $\nu$, $m$ and $\delta_2$) and a unique local classical  solution
$(\eta, u,q)\in C^0(\overline{I_T},{H}^{18} )\times \mathcal{U}_{T}\times C^0(\overline{I_T},\underline{H}^{17})$
to the initial value problem \eqref{01dsaf16asdfasf}--\eqref{01dsaf16asdfasfsaf}, satisfying
$0<\inf_{(y,t)\in \mathbb{R}^3\times \overline{I_T}} \det(\nabla \eta+I)$ and $\sup_{t\in \overline{I_T}}\|\nabla \eta\|_{12}\leqslant 2\delta_2$\footnote{Here the uniqueness means that if there is another solution
$( \tilde{\eta}, \tilde{u},\tilde{q})\in C^0(\overline{I_T},{H}^{18} ) \times \mathcal{U}_{T}\times C^0(\overline{I_T},\underline{H}^{17})$
 satisfying $0<\inf_{(y,t)\in \mathbb{R}^3\times \overline{I_T}} \det(\nabla \tilde{\eta}+I)$, then
 $(\tilde{\eta},\tilde{u},\tilde{q})=(\eta,u,q)$ by virtue of the smallness condition
 ``$\sup_{t\in \overline{I_T}}\|\nabla \eta\|_{12}\leqslant 2\delta_2$''.}.
\end{pro}
\begin{pf}
Please refer to  \cite[Proposition 7.3]{JFJSAAD}.
\hfill $\Box$
\end{pf}
\begin{pro}\label{pro:0401nasfxdxx}
There is a positive constant $\delta_3 $, such that for any $\varphi\in H^{18}$ satisfying
$\|\nabla \varphi\|_{2}\leqslant \delta_3$, we have (after possibly being redefined on a set of {zero measure})
$\det(\nabla \varphi+I)>1/2$ and
\begin{align}
 \psi : \mathbb{R}^3\to \mathbb{R}^3 \mbox{ is a }C^{16}\mbox{ homeomorphism mapping},
\end{align}
where $\psi:=\varphi+y$.
\end{pro}
\begin{pf}
 Please refer to \cite[Lemma 4.2]{JFJSOMITIN} for a detailed proof.
\hfill $\Box$
\end{pf}

With the \emph{a priori} estimate \eqref{20202110261936} (under the conditions \eqref{aprpiosesnewxxxx} and \eqref{aprpiose1snewxxxxz}
with $\delta\leqslant \delta_1$), and Propositions \ref{pro:0401nxdxx}, \ref{pro:0401nasfxdxx} in hand, we can easily establish
Theorem \ref{201904301948}. We briefly give the proof below.

Let $m$ and $(\eta^0,u^0)\in (\underline{H}^{18}_{1}\cap H^{18}_*) \times \underline{H}^{17}$ satisfy
\begin{align}
\max\{K ,K^2\}/m\leqslant\min\{\delta_1,  \delta_2/c_{\omega,13}, \delta_3/c_{\omega,13
}\}=:c_3\leqslant 1,
\label{201912261425}
\end{align}
where $K$ is defined by \eqref{201911262060}, and the above two constants $c_{\omega,13}$ come  from  \eqref{2020211013200910} with $i=13$.
Then we see that $\eta^0$ satisfies \eqref{201912261426} by \eqref{2020211013200910} and \eqref{201912261425}. Hence,
by virtue of Proposition \ref{pro:0401nxdxx}, there exists a unique local solution $(\eta, u,q)$
of \eqref{01dsaf16asdfasf}--\eqref{01dsaf16asdfasfsaf} with  {the} maximal existence time $T^{\max}$, satisfying
\begin{itemize}
  \item for any $T\in  I_{T^{\max}}$,
the solution $(\eta,u,q)$ belongs to $C^0(\overline{I_T},\underline{H}^{18}_1)\times \underline{\mathcal{U}}_T\times C^0( \overline{I_T},\underline{H}^{17}) $ and $$\sup_{t\in \overline{I_T}} \|\nabla \eta\|_{12}\leqslant 2\delta_2;$$
  \item $\limsup_{t\to T^{\max} }\|\nabla \eta( t)\|_{12} > \delta_2$ or $\limsup_{t\to T^{\max} }\|(\nabla \eta,u)( t)\|_{17}=\infty$, if $T^{\max}<\infty$.
\end{itemize}

Let
\begin{align}
&E(t):=  \|\eta(t)\|_{18}^2 +  m^{-2/3} \|  u (t) \|_{17}^2 +
 \| m\partial_\omega \eta (t)\|_{16}^2,
\nonumber\\
&T^{*}=\sup\left\{ T \in I_{T^{\max}}~\left|~ E(t)\leqslant K^2 \mbox{ for any }t\leqslant T\right.\right\}.\nonumber
\end{align}
Recalling the definition of $K$ and  the condition $c_1\geqslant 4$, we easily see that the definition of $T^*$ makes sense and $T^*>0$.
In addition,
by \eqref{2020211013200910}, we have $\| \nabla \eta\|_{12} \leqslant \delta_3$ for all $t\in I_{T^*}$,
then $\eta(t)\in H^{18}_*$ for all $t\in I_{T^*}$ by Proposition \ref{pro:0401nasfxdxx}.
Thus, to obtain the existence of a global solution, it suffices to verify $T^*=\infty$. Now, we show this by contradiction.

Assume $T^*<\infty$. Keeping in mind that $T^{\max}$ denotes the maximal existence time and $K/m\leqslant \delta_2/c_{\omega,13}$ by virtue of
\eqref{201912261425}, we apply Proposition \ref{pro:0401nxdxx} to find that $T^{\max}>T^*$  and
\begin{equation}
\label{201911262202}
 E(T^*)=K^2.
\end{equation}
Since $m^{-1}\max\{K,K^2\}\leqslant  \delta_1$ and $\sup_{0\leqslant t\leqslant T^*} {E} (t) \leqslant K^2$,
 we can still show that the solution $(\eta,u )$ enjoys the stability estimate \eqref{20202110261936} with $T^*$ in place of $T$
 by the regularity of $(\eta,u,q)$. More precisely, we have $\sup_{0\leqslant t\leqslant T^*}  {E} (t) \leqslant K^2/4$,
which contradicts with \eqref{201911262202}. Hence, $ T^*=\infty$, and thus $T^{\max}=\infty$.

The uniqueness of the global solutions is obvious due to the uniqueness of the local solutions in Proposition \ref{pro:0401nxdxx}
 and the fact $\sup_{t\geqslant 0}\|\nabla \eta\|_{12}\leqslant 2 \delta_2$. Finally, it is obvious that the global solution $(\eta,u,q)$ enjoys the estimates \eqref{20190safd5041053}--\eqref{20202111005120114}, \eqref{200s12}--\eqref{20200sasf12} and \eqref{2017020614181721} by recalling the derivation of  \emph{a priori} energy estimates for $(\eta,u)$. This completes the proof of Theorem \ref{201904301948}.

\section{Proof of Theorem \ref{201912041028}}\label{2020011sdf92326}

This section is devoted to the proof of Theorem \ref{201912041028}.
 Let $(\eta^0,u^0)$ satisfy all the assumptions in Theorem \ref{201904301948} and $(\eta,u,q)$ be the solution constructed by Theorem \ref{201904301948}.
 Exploiting \eqref{201909281832}, \eqref{20202111005120114} and \eqref{2020211013200910}, we have
\begin{align}
\| \eta^0\|_{15}\lesssim 1\mbox{ and }\|\eta^0\|_{13}\sqrt{\mathcal{E}_{H}^0} \lesssim 1.
\label{20211111311}
\end{align}

 By the regularity theory of the Stokes problem, there exists a unique solution $(\eta^{\mm{r}},u^{\mm{r}},Q_1,Q_2)$ satisfying
\begin{align}
   \begin{cases}
   -\Delta \eta^{\mm{r}} +\nabla Q_1=0, \\
   \mm{div} \eta^{\mm{r}}=-\mm{div}\eta^0,\\
     (\eta^{\mm{r}})_{\mathbb{T}^3}  = 0
   \end{cases}\mbox{and } \begin{cases}
-\Delta u^{\mm{r}} +\nabla Q_2=0, \\
\mm{div} u^{\mm{r}}=\mm{div}_{\tilde{\mathcal{A}}^0}u^0,\\
({u}^{\mm{r}})_{\mathbb{T}^3}  = 0, \nonumber
\end{cases}
\end{align}
where $\tilde{\mathcal{A}}^0:={\mathcal{A}}^0-I$. Moreover, $(\eta^{\mm{r}},u^{\mm{r}})\in   \underline{H}^{18} \times \underline{H}^{17} $ satisfies  \eqref{202af011032122} and \eqref{202011032123} by making use of the classical regularity theory of Stokes equations, the integral by parts, Poinc\'are's inequality, Young's inequality, \eqref{201909261909}, the interpolation inequality \eqref{202108261405},  \eqref{200s12} with $t=0$, \eqref{2020saf05sda08}, \eqref{202005sda08}, \eqref{202108271645}, \eqref{20211111311} and the following identity
 \begin{align*}
\mm{div}\eta^0=\mm{div} \left(\begin{array}{c}
         - \eta_1^0(\partial_2\eta_2^0+\partial_3\eta_3^0 )+
           \eta_1^0(\partial_2\eta_3^0 \partial_3\eta_2^0
- \partial_2\eta_2^0\partial_3\eta_3^0) \\
      \eta_1^0\partial_1\eta_2^0  - \eta_2^0\partial_3\eta_3^0+
        \eta_1^0( \partial_1\eta_2^0 \partial_3\eta_3^0
-\partial_1\eta_3^0\partial_3\eta_2^0)\\
          \eta_1^0\partial_1\eta_3^0
+\eta_2^0\partial_2\eta_3^0+\eta_1^0(\partial_1\eta_3^0 \partial_2\eta_2^0
-\partial_1\eta_2^0\partial_2\eta_3^0)
        \end{array}\right).
 \end{align*}
 Moreover, using the interpolation inequality \eqref{202108261405}, we have
\begin{align}
& \begin{cases}
\|\eta^{\mm{r}}\|_{k+1}\lesssim  \|\eta^0\|_3\| \eta^0\|_{k+1 }\mbox{ for }0\leqslant k\leqslant 17,&\\
 \|\partial_\omega \eta^{\mm{r}}\|_k\lesssim \|  \eta^0\|_{3} \|\partial_\omega \eta^0\|_{k }   +
\begin{cases}
 \| \eta^0\|_7\| \partial_\omega \eta^0\|_{k}&\mbox{for }0\leqslant k\leqslant 2 ;\\
\| \eta^0\|_k\| \partial_\omega \eta^0\|_3&\mbox{for } 3\leqslant k\leqslant 17,
 \end{cases}&  \\
  \| u^{\mm{r}} \|_{k}\lesssim \|  \eta^0\|_3\|u^0\|_k +
\begin{cases}
  \sqrt{\|\partial_\omega \eta^0\|_8\|  u^0\|_6\| \partial_\omega\eta^0\|_k \|  u^0\|_{k}}&\mbox{for }0\leqslant k\leqslant 2 ;\\
\|  \eta^0\|_{k }\|  u^0\|_{3 }&\mbox{for }3\leqslant k\leqslant 17.
 \end{cases}&
 \end{cases}
  \label{2032123}
\end{align}

Let $\tilde{\eta}^0=\eta^0+\eta^{\mm{r}}$ and $\tilde{u}^0=u^0+u^{\mm{r}}$. Thus, it is easy to see that $(\tilde{\eta}^0,\tilde{u}^0)$ belongs to $ \underline{H}^{18}_\sigma\times \underline{H}^{17}_\sigma$. Therefore, there exists a unique global solution $(\eta^{\mm{L}},u^{\mm{L}},q^{\mm{L}} )\in C^0(\mathbb{R}_0^+,\underline{H}^{18}_\sigma)\times \underline{\mathcal{U}}_\infty\times C^0(\mathbb{R}_0^+,\underline{H}^{17})$ to the linearized problem \eqref{202001070914} with the initial condition $(\eta^{\mm{L}} ,u^{\mm{L}})|_{t=0}=(\tilde{\eta}^0,\tilde{u}^0 )$.

Similarly to \eqref{2019112safdsafsaf6201sdaf2} and \eqref{20191201sdaf2}, we  easily see that the solution $(\eta^{\mm{L}}, u^{\mm{L}})$ of the linearized problem \eqref{202001070914}  enjoys the estimates \eqref{2020052010sdfa33} and \eqref{202033}.
Moreover, by \eqref{20211111311} and \eqref{2032123}, we can further derive form \eqref{2020052010sdfa33} and \eqref{202033} that
\begin{align}
  \mathcal{E}_{j}^{\mm{L}} +
\int_0^t\mathcal{D}_{j}^{\mm{L}} \mm{d}\tau\lesssim \mathcal{E}^0_{j}\mbox{ for any }1\leqslant j\leqslant 15 \label{033}
\end{align}
and
\begin{align}
\sum_{i=0}^{3} \left((1+m^{-2})^i \langle t\rangle^{ (3-i)} {\mathcal{E}}_{4i}^{\mm{L}}+ (1+m^{-2})^i\int_0^t\langle t\rangle^{ (3-i)}\mathcal{D}_{4i}^{\mm{L}}  \mm{d}\tau\right)\lesssim \Xi , \label{2023}
\end{align}
where $\Xi$ is defined by \eqref{20202111110854}.

Let $(\eta^{\mm{d}}, u^{\mm{d}})= (\eta-\eta^{\mm{L}}, u-u^{\mm{L}})$, then the error function $(\eta^{\mm{d}}, u^{\mm{d}} )$ satisfies
\begin{equation}\label{01dsaf16safafasdfasfxx}
\begin{cases}
\eta_t^{\mm{d}}=u^{\mm{d}}, \\[1mm]
u_t^{\mm{d}} -\nu\Delta    u^{\mm{d}}- m^2
\partial_\omega^2\eta^{\mm{d}}=\mathfrak{N}, \\[1mm]
\div u^{\mm{d}}=-\mathrm{div}_{\tilde{\mathcal{A}}} {u},\\
(\eta^{\mm{d}},u^{\mm{d}})|_{t=0}=-(\eta^{\mm{r}},u^{\mm{r}}).
\end{cases}
\end{equation}
It is easy to see from \eqref{01dsaf16safafasdfasfxx} that $(u^{\mm{d}})_{\mathbb{T}^3}=(\eta^{\mm{d}})_{\mathbb{T}^3}=0$, $\mm{div}\eta^{\mm{d}}=\mm{div}\eta$  for any $t> 0$, since $(\eta^{\mm{r}})_{\mathbb{T}^3}  = ({u}^{\mm{r}})_{\mathbb{T}^3}=0$ and $\mm{div}\eta^{\mm{L}}= 0$.

Recalling that $(\eta,u,q)$ is constructed by Theorem \ref{201904301948}, the solution $(\eta,u,q)$ satisfies   all the estimates
in \eqref{20202111005120114},  \eqref{200s12}--\eqref{20200sasf12} and \eqref{2017020614181721}. Hence, we can follow the  {same} arguments as in the  proof of
Lemmas \ref{qwe201612132242nn}--\ref{qwe201612132242nxsfssdfsxx}
with some modifications to deduce from \eqref{01dsaf16safafasdfasfxx} that
\begin{align}
  \frac{\mm{d}}{\mm{d}t}  \| \nabla  (u^{\mm{d}},m \partial_\omega \eta^{\mm{d}})\|_i^2 +
 \nu\|  u^{\mm{d}} \|_{i+1}^2 \lesssim \varphi_i/m,\label{201702061safasaf4181721asfdx}
 \end{align}
  and
  \begin{align}
&  \frac{\mm{d}}{\mm{d}t}\bigg( \frac{\nu}{2}\|\nabla^{i+1} \eta^{\mm{d}}\|_0^2+ \sum_{|\alpha|=i} \int \partial^\alpha \eta^{\mm{d}} \cdot \partial^\alpha u^{\mm{d}}\mm{d}y  \bigg)+\|m \partial_\omega \eta^{\mm{d}} \|_i^2 \leqslant   \| u^{\mm{d}} \|_i^2  +c \psi_i/m, \label{Lem:03sfdaxx}\end{align}
where
  \begin{align} \varphi_i :=&m\| (u,u^{\mm{L}})\|_{i+1} ( \|\partial_\omega \eta\|_6 \|u\|_{i+1} +\| \partial_\omega \eta\|_{i+4}\|u\|_3)\nonumber \\
 &+ m\begin{cases}
  \| \partial_\omega  \eta\|_{7 } \|(u,u^{\mm{L}})\|_4 ( \|  u\|_3^2+ m^2 \| \partial_\omega   \eta \|_3 \|  \partial_\omega \eta \|_6 ) &\mbox{for }0\leqslant i\leqslant 3;  \\
 (\|\partial_\omega \eta \|_6\| (u,u^{\mm{L}})\|_{i+1} +\|\partial_\omega \eta\|_{i+4 } \|(u,u^{\mm{L}})\|_3)    &\\
( \|   u\|_3 \| u\|_{i}
 + m^2 (\| \partial_\omega \eta \|_6 \| \partial_\omega \eta \|_{i } + \| \partial_\omega  \eta \|_4\|
\partial_\omega \eta \|_{i+2}
))  &\\
+\|  \partial_\omega \eta\|_{i+3} \| (u,u^{\mm{L}})\|_{i+1}(\|  u\|_3^2+ m^2 \| \partial_\omega   \eta \|_3 \|  \partial_\omega \eta \|_6)    &\mbox{for }4\leqslant i\leqslant 12
  \end{cases}  \label{2022111091024}
  \end{align}
and
 \begin{align} \psi_i:= &m \| (\eta, \eta^{\mm{L}})\|_{i+1}  (\|\partial_\omega \eta \|_6 \|u\|_{i+1} + \|\partial_\omega \eta\|_{i+4}\|u\|_3 )\nonumber \\
 &+  m\|(\eta, \eta^{\mm{L}})\|_{i+1} \begin{cases}
\| \partial_\omega \eta\|_6(\| u\|_3^2
  +m^2    \| \partial_\omega   \eta \|_3 \|  \partial_\omega \eta \|_6 ) &\mbox{for }0\leqslant i\leqslant 3;  \\
 \|\partial_\omega \eta \|_6 ( \|  u\|_3 \|   u\|_{i}+m^2 \| \partial_\omega \eta \|_6\| \partial_\omega \eta \|_{i}) &  \\
+ \|\partial_\omega  \eta\|_{i+3}  ( \| u\|_3^2  + m^2\| \partial_\omega \eta \|_4  \| \partial_\omega \eta \|_6 )&\mbox{for }4\leqslant i\leqslant 12 .  \label{2024}
  \end{cases}
\end{align}

Similarly to \eqref{20191sadf1262012}, we further derive from \eqref{201702061safasaf4181721asfdx} and \eqref{Lem:03sfdaxx} that
\begin{align}
&\frac{\mm{d}}{\mm{d}t}    \tilde{\mathcal{E}}_{i}^{\mm{d}}+
c \mathcal{D}_{i}^{\mm{d}}    \lesssim m^{-1}(\varphi_i+ \psi_i) \mbox{ for }0\leqslant i\leqslant 12,
\label{2016sfd2012}
\end{align}
where $
\tilde{\mathcal{E}}_{i}^{\mm{d}} $ is defined as $\tilde{\mathcal{E}}_{i} $ in \eqref{202211052115} with $(\eta^{\mm{d}},u^{\mm{d}})$ in place of $(\eta ,u )$,
and $\tilde{\mathcal{E}}_{i}^{\mm{d}} $ satisfies
\begin{align}
 \mathcal{E}_{i}^{\mm{d}}\lesssim   \tilde{\mathcal{E}}_{i}^{\mm{d}}\lesssim  \mathcal{E}_{i}^{\mm{d}}. \label{20200f6141102}
\end{align}
In addition, by \eqref{2053}, \eqref{202211111}, \eqref{201041053} and \eqref{033}, it is easy to see that
\begin{align}
&\int_0^t   (\varphi_{12}+ \psi_{12}) \mm{d}\tau \nonumber \\
&\lesssim  m\int_0^t  ( \| (u,u^{\mm{L}})\|_{13} ( \|\partial_\omega \eta\|_6 \|u\|_{13} +\| \partial_\omega \eta\|_{16}\|u\|_3+\|  \partial_\omega \eta\|_{15}  ( \|  u\|_3^2\nonumber\\
& \quad+ m^2 \| \partial_\omega   \eta \|_3 \|  \partial_\omega \eta \|_6)) +(\|\partial_\omega \eta \|_6\| (u,u^{\mm{L}})\|_{13} +\|\partial_\omega \eta\|_{16 } \|(u,u^{\mm{L}})\|_3)\nonumber      \\
&\quad  ( \|   u\|_3 \| u\|_{12}
 + m^2 (\| \partial_\omega \eta \|_6 \| \partial_\omega \eta \|_{12} + \| \partial_\omega  \eta \|_4 \|
\partial_\omega \eta \|_{14}
))
 \nonumber\\
 &\quad + \| (\eta, \eta^{\mm{L}})\|_{13} (  \|\partial_\omega \eta \|_6 \|u\|_{13} + \|\partial_\omega \eta\|_{16}\|u\|_3  +\|\partial_\omega \eta \|_6   ( \|  u\|_3 \|   u\|_{12}\nonumber \\
 &\quad +m^2 \| \partial_\omega \eta \|_6\| \partial_\omega \eta \|_{12}) +\|\partial_\omega  \eta\|_{15} (  \| u\|_3^2 + m^2\| \partial_\omega \eta \|_4  \| \partial_\omega \eta \|_6 )  ))\mm{d}\tau\nonumber \\
&\lesssim \int_0^t\left(  \sqrt{\mathcal{E}_{16}+\mathcal{E}_{16}^{\mm{L}}
} \left(1+\sqrt{\mathcal{E}_{16}} \right)  (\mathcal{D}_{12} +
 \mathcal{D}_{12}^{\mm{L}}  ) + \sqrt{ \mathcal{E}_3\mathcal{E}_{16} (\mathcal{E}_{12} + \mathcal{E}_{12}^{\mm{L}} )} \right)\mm{d}\tau\nonumber \\
 & \lesssim{\mathcal{E}^{0}_{17}
e^{c_2 \vartheta}} \left( \sqrt{\mathcal{E}_{12}^0}+ {\mathcal{E}_{12}^0} \right)  +\sqrt{\Xi\mathcal{E}_{12}^0 \mathcal{E}^{0}_{17}
e^{c_2 \vartheta}}. \label{20221110111029}
\end{align}
Thus, integrating \eqref{2016sfd2012} over $(0,t)$, and then making use of \eqref{20211111311}, \eqref{2032123}, \eqref{20200f6141102} and \eqref{20221110111029}, we easily get \eqref{2020sfa05201033}.

Finally we derive \eqref{201905041053xxx}. Similarly to  \eqref{20191201sdaf2}, we  derive the following inequality from  \eqref{2020211013200910}, \eqref{2032123}  \eqref{2016sfd2012} and \eqref{20200f6141102}:
\begin{align}
 &\sum_{i=0}^{3} \left((1+m^{-2})^i \langle t\rangle^{ (3-i)} {\mathcal{E}}_{4i}^{\mm{d}}+ (1+m^{-2})^i \int_0^t\langle \tau\rangle^{ (3-i)}\mathcal{D}_{4i}^{\mm{d}}\mm{d}\tau  \right)\nonumber \\
 &\lesssim
 m^{-1}\left(\mathcal{E}_{15}^0 \Xi +\sum_{i=0}^{3} \int_0^t  (1+m^{-2})^i\langle \tau \rangle^{3-i}(\varphi_i+ \psi_i))\mm{d}\tau \right).
\label{20sdfa112}
\end{align}

Noting that $(\eta,u)$ satisfies \eqref{2x62012}, thus making use of  the interpolation inequality \eqref{202004221saffad412}, \eqref{20190safd5041053}, \eqref{2053}, \eqref{033} and \eqref{2023}, we have
\begin{align}
\nonumber
&\int_0^t  \langle \tau\rangle^{3 }(\varphi_0+ \psi_0) \mm{d}\tau\\
&\lesssim m\int_0^t \langle \tau\rangle^{3 }(\| (u,u^{\mm{L}})\|_{1} ( \|\partial_\omega \eta\|_6 \|u\|_{1} +\| \partial_\omega \eta\|_{4}\|u\|_3)  +
 \|\partial_\omega  \eta\|_{7}\| (u,u^{\mm{L}})\|_4( \|  u \|_3^2\nonumber \\
&\quad+
   m^2  \|\partial_\omega \eta\|_3\|\partial_\omega \eta\|_6 ) +  \| (\eta, \eta^{\mm{L}})\|_{1}  (\|\partial_\omega \eta \|_6 \|u\|_{1} + \|\partial_\omega \eta\|_{4}\|u\|_3 )\nonumber \\
 &\quad+
 \|(\eta, \eta^{\mm{L}})\|_1\| \partial_\omega \eta\|_6(  \| u\|_3^2 +m^2  ( \| \partial_\omega   \eta \|_2 \|  \partial_\omega \eta \|_3+ \| \partial_\omega \eta \|_3\| \partial_\omega \eta \|_6 ))) \mm{d}\tau \nonumber\\
&\lesssim\int_0^t \langle \tau\rangle^{3 }   ( 1+  \sqrt{\mathcal{E}_{12}})
  \sqrt{(\mathcal{E}_{0} + \mathcal{E}_{0}^{\mm{L}})( \mathcal{D}_{0}+\mathcal{D}_{0}^{\mm{L}} )(\mathcal{D}_{12}+\mathcal{D}_{12}^{\mm{L}} )} \mm{d}\tau \nonumber \\ &\lesssim\Xi \left(\sqrt{\mathcal{E}_{12}^0}+{\mathcal{E}_{12}^0}\right).\nonumber
\end{align}
Similarly,
\begin{align}
&\int_0^t  \langle \tau\rangle^{2 }(\varphi_4+ \psi_4) \mm{d}\tau\nonumber\\
&\lesssim m \int_0^t \langle \tau\rangle^{2 }(\| (u,u^{\mm{L}})\|_{5} ( \|\partial_\omega \eta\|_6 \|u\|_{5} +\| \partial_\omega \eta\|_{8}\|u\|_3\nonumber\\
& \quad + \|\partial_\omega \eta\|_{8} (  \|  u \|_{4}^2+ m^2 \| \partial_\omega\eta \|_6 \|  \partial_\omega \eta \|_{4} ))  + \| (\eta, \eta^{\mm{L}})\|_{5} ( \|\partial_\omega \eta \|_6 \|u\|_{5} \nonumber\\
 &\quad + \|\partial_\omega \eta\|_{8}\|u\|_3 +\|\partial_\omega  \eta\|_{7} (
      \| u\|_4^2  + m^2\| \partial_\omega \eta \|_4  \| \partial_\omega \eta \|_6 )  )) \mm{d}\tau\nonumber \\
&\lesssim \int_0^t \langle \tau\rangle^{2 } ( 1+  \sqrt{\mathcal{E}_{12}})
  \sqrt{\mathcal{D}_{12}(\mathcal{E}_{4} + \mathcal{E}_{4}^{\mm{L}})( \mathcal{D}_{4}+\mathcal{D}_{4}^{\mm{L}} )}\mm{d}\tau\nonumber
\\ &\lesssim\Xi\left( \sqrt{\mathcal{E}_{12}^0}+\mathcal{E}_{12}^0\right) \nonumber
\end{align}
and
\begin{align}
&\int_0^t  \langle \tau\rangle (\varphi_8+ \psi_8) \mm{d}\tau\nonumber\\
&\lesssim m\int_0^t\langle \tau\rangle (\| (u,u^{\mm{L}})\|_{9} ( \|\partial_\omega \eta\|_6 \|u\|_{9} +\| \partial_\omega \eta\|_{12}\|u\|_3+\|  \partial_\omega \eta\|_{11}  (  \|  u\|_3^2 \nonumber\\
& \quad + m^2 \| \partial_\omega   \eta \|_3 \|  \partial_\omega \eta \|_6))+( \|   u\|_3 \| u\|_8 + m^2 (\| \partial_\omega \eta \|_6 \| \partial_\omega \eta \|_8\nonumber \\
 &\quad
+\| \partial_\omega  \eta \|_4 \|
\partial_\omega \eta \|_{10}
)) (\|\partial_\omega \eta \|_6\| (u,u^{\mm{L}})\|_{9} +\|\partial_\omega \eta\|_{12 } \|(u,u^{\mm{L}})\|_3)  \nonumber\\
 &\quad + \| (\eta, \eta^{\mm{L}})\|_{9} (  \|\partial_\omega \eta \|_6 \|u\|_{9} + \|\partial_\omega \eta\|_{12}\|u\|_3   +\|\partial_\omega  \eta\|_{11} (  \| u\|_3^2 \nonumber \\
 & \quad+ m^2\| \partial_\omega \eta \|_4  \| \partial_\omega \eta \|_6 )  +\|\partial_\omega \eta \|_6   ( \|  u\|_3 \|   u\|_{8}+m^2 \| \partial_\omega \eta \|_6\| \partial_\omega \eta \|_{8}))) \mm{d}\tau\nonumber\\
&\lesssim \int_0^t \langle \tau\rangle ( 1+  \sqrt{\mathcal{E}_{12}})
  \sqrt{\mathcal{D}_{12}(\mathcal{E}_{8} + \mathcal{E}_{8}^{\mm{L}})( \mathcal{D}_{8}+\mathcal{D}_{8}^{\mm{L}} )} \mm{d}\tau \nonumber \\
&\lesssim\Xi  \left(
\sqrt{\mathcal{E}_{12}^0}+\mathcal{E}_{12}^0
\right). \nonumber
\end{align}

Putting the above three estimates and \eqref{20221110111029} together yields
\begin{align}
\nonumber
&\sum_{i=0}^{3} \int_0^t  (1+m^{-2})^i\langle \tau\rangle^{3-i}(\varphi_i+ \psi_i))\mm{d}\tau \lesssim \Phi ,\nonumber
\end{align}
where $\Phi$ is defined by \eqref{2020211111342}.
Inserting the above estimate into \eqref{20sdfa112} yields \eqref{201905041053xxx}. This completes the proof of Theorem \ref{201912041028}.

%x\section{Proof of Proposition \ref{20206141015}}

\vspace{4mm} \noindent\textbf{Acknowledgements.}
% {The authors would like to thank the anonymous referees for valuable suggestions and comments, which improve the presentation of this paper.}
The research of Fei Jiang was supported by NSFC (Grant Nos.  12022102) and the Natural Science Foundation of Fujian Province of China (2020J02013), and the research of Song Jiang by National Key R\&D Program (2020YFA0712200), National Key Project (GJXM92579), and
NSFC (Grant No. 11631008), the Sino-German Science Center (Grant No. GZ 1465) and the ISF--NSFC joint research program (Grant No. 11761141008).
\renewcommand\refname{References}
\renewenvironment{thebibliography}[1]{ %
\section*{\refname}
\list{{\arabic{enumi}}}{\def\makelabel##1{\hss{##1}}\topsep=0mm
\parsep=0mm
\partopsep=0mm\itemsep=0mm
\labelsep=1ex\itemindent=0mm
\settowidth\labelwidth{\small[#1]}%
\leftmargin\labelwidth \advance\leftmargin\labelsep
\advance\leftmargin -\itemindent
\usecounter{enumi}}\small
\def\newblock{\ }
\sloppy\clubpenalty4000\widowpenalty4000
\sfcode`\.=1000\relax}{\endlist}
\bibliographystyle{model1b-num-names}
 
\end{document}